\documentclass[11pt]{article}
\usepackage{amsgen,amsmath,amstext,amsbsy,amsopn, amssymb,amsfonts}
\usepackage{pstricks}
\usepackage{epsf}
\usepackage{float}
\usepackage{epsfig}

\usepackage{fullpage}

\usepackage{amsfonts}
\usepackage{amssymb}

\newcommand{\iref}[1]{\eqref{#1}}

\def \be{\begin{equation}}
\def \ee{\end{equation}}
\newcommand{\beq}{\begin{eqnarray}}
\newcommand{\eeq}{\end{eqnarray}}
\newcommand{\bea}{\begin{array}{c}}
\newcommand{\eea}{\end{array}}
\newcommand{\bi}{\begin{itemize}}
\newcommand{\ei}{\end{itemize}}

\newtheorem{lem}{Lemma}

\newtheorem{theo}{Theorem}

\def\N{\mathbb N}
\def\R{\mathbb R}

\def\C{\mathbb C}

\def\vp{\varphi}
\def\t{\tilde}
\def\e{\varepsilon}
\def\eps{\varepsilon}

\def\Chi{\raise .3ex \hbox{\large $\chi$}} 
\def\vp{\varphi}
\def\s{\sigma}
\def\d{\delta}

\def\b{\beta}

\def\p{\partial}

\def\ii{\infty}

\def\BP{\mathcal{P}}

\def\BE{\mathcal{E}}

\def\BC{\mathcal{C}}
\def\BO{\mathcal{O}}
\def\BI{\mathcal{I}}
\def\ra{\to}

\def\b{\ \ \ \ $\Box$}
\def\ds{\displaystyle}
\def\bs{\boldsymbol}

\title{The KdV/KP-I limit of the Nonlinear Schr\"odinger equation}
\author{D. Chiron\footnote{Laboratoire J.A. DIEUDONNE, Universit{\'e} 
de Nice - Sophia Antipolis, Parc Valrose, 06108 Nice Cedex 02, France. 
{\it e-mail:} \texttt{chiron@unice.fr}} \,\, \& \,\, F. Rousset.\footnote{IRMAR, Universit\'e de Rennes 1, 
Campus de Beaulieu, 35042 Rennes Cedex, France. {\it e-mail:} \texttt{ frederic.rousset@univ-rennes1.fr}}}
\date{}

\begin{document}

\maketitle
\begin{abstract}
We justify rigorously the convergence of the amplitude of solutions of 
Nonlinear-Schr\"odinger type  Equations with non zero limit at infinity  to 
an asymptotic regime governed by the 
Korteweg-de Vries equation in dimension $1$ and the Kadomtsev-Petviashvili I 
equation in dimensions $2$ and more. We  get two types of results.  In the 
one-dimensional case, we   prove directly  by energy bounds that there is no vortex formation
 for  the global solution of the NLS 
equation in the energy space and deduce from this the convergence  towards the unique solution in the energy 
space of the KdV equation. In arbitrary dimensions, we use an hydrodynamic reformulation 
of NLS and recast the problem as a singular limit for an hyperbolic system. We thus 
prove that smooth $H^s$ solutions  exist on a time interval independent of the small 
parameter. We then pass to the limit by a compactness argument and obtain the KdV/KP-I equation.

\end{abstract}

\section{Introduction}

\ \indent We consider the $n$-dimensional nonlinear Schr\"odinger equation
\be
\tag{NLS}
i \frac{\p \Psi}{\p \tau} + \frac{1}{2} \Delta_z \Psi = \Psi f(|\Psi|^2) 
\quad \quad \quad \quad \Psi = \Psi(\tau, z) : \R_+ \times \R^n \to \C.
\ee
This equation is used as a model in nonlinear Optics (see for instance 
\cite{KL}) and in superfluidity and Bose-Einstein condensation 
(see, {\it e.g.} \cite{RB}, \cite{GP}, \cite{G}).

We assume that, for some $\rho_0>0$, $f(\rho_0^2)=0$, so that $\Psi \equiv \rho_0$ 
is a particular solution of  (NLS). We are interested in solutions $\Psi$ of (NLS) 
such that $|\Psi| \simeq \rho_0$. In the sequel, we take $\rho_0 =1$, 
the general case follows changing $\Psi$ for $\t{\Psi} \equiv \rho_0^{-1} \Psi$ 
and $\t{f}(R) \equiv f(\rho_0^2 R )$. Then, from now on, we consider smooth nonlinearities
 $f \in \mathcal{C}^\infty(\mathbb{R}, \mathbb{R})$ such that 
\be
\label{hyp} f(1)=0, \quad \quad f'(1)>0 
\ee
and will be interested in situations where  $| \Psi| \simeq 1$. 
 Note that  this means thanks to \eqref{hyp} that we shall study the equation in a defocusing regime.
A typical example of nonlinearity  is 
simply $f(R) = R-1$ for which  (NLS)  is termed the Gross-Pitaevskii equation. Equation 
 (NLS) is  an  Hamiltonian  flow associated to the Ginzburg-Landau type energy (when it makes sense)
$$ \BE(\Psi) \equiv 
\frac{1}{2} \int_{\R^n} |\nabla_{z} \Psi|^2 + F \big( |\Psi|^2 \big) \ dz, $$
where $F(R) \equiv 2 \ds{\int_1^R f(r) \ dr} $.

\subsection{KdV and KP-I asymptotic regimes for NLS}

\ \indent In a suitable scaling corresponding to $| \Psi| \simeq 1$, the dynamics 
for the amplitude of $\Psi$ converges, in dimension $n=1$, to the Korteweg-de Vries 
equation
\be
\tag{KdV}
2 \p_t v + k \, v \p_x v - \frac{1}{4c^2}\, \p_{xxx} v = 0,
\ee
and in dimensions $n\geq 2 $ to the Kadomtsev-Petviashvili - I equation
\be
\tag{KP-I}
\p_x \Big( 2 \p_t v + k v \p_x v - \frac{1}{4c^2} \, \p ^3_x v \Big) + \Delta_{\perp} v = 0
\ee
where $v=v(t,X) \in \R$, $X=(x,x_\perp)\in \R \times \R^{n-1}$.
 The coefficients $c$ and $k$ are related to the nonlinearity $f$ by
\be
\label{defcoeff}
c \equiv \sqrt{f'(1)} \quad \quad \quad \quad {\rm and} \quad \quad \quad \quad 
k \equiv 6 + \frac{2}{c^2} f''(1).
\ee
Note that the KP-I equation   reduces to the KdV equation if $v$ does not depend on $x_{\perp}$.

The formal derivation of this regime is as follows. First, we consider 
a small parameter $\e$, and rescale time and space according to
\begin{equation}
\label{scaling}
t = c \e^3 \tau, \quad  X_{1} = x = \e (z_1 - c \tau), \quad X_j = \e^2 z_j, 
\quad j \in \{ 2, ... , n\},  \quad  \Psi(\tau,z) = \psi^\e(t, X).
\end{equation}
The nonlinear Schr\"odinger equation  for  $\psi^\e$ reads now
\be
\label{NLSd}
 i c \e^3 \frac{\p \psi^\e}{\p t} 
- i c \e \p_x \psi^\e 
+ \frac{\e^2}{2} \p_x^2 \psi^\e 
+ \frac{\e^4}{2} \Delta_{\perp} \psi^\e
= \psi^\e f(|\psi^\e|^2), \quad \quad \quad X= (x,x_\perp) \in \R \times \R^{n-1}.
\ee
We shall use  the following ansatz for $\psi^\e$
\be
\label{ansatz}
\psi^\e(t, X) = \big( 1+ \e^2 A^\e (t,X)\big) \exp\big( i \e \vp^\e(t,X) \big)
\ee
where the amplitude $A^\eps \in \mathbb{R}$ is assumed to be of order $1$ and
 the real phase $\varphi^\eps \in \mathbb{R}$ is also assumed to be of order $1$.
 The ansatz \eqref{scaling}, \eqref{ansatz} mean  that we study a  weak amplitude  wave propagating to the right  in a long wave regime and
 that this wave is slowly modulated in the transverse direction
  thanks to \eqref{scaling}.
 Note that the  occurence  of the KdV or KP equation as enveloppe equations  in such regimes  
 is expected. We refer for example to \cite{AL} and references therein
  for the derivation of    these equations
  from the water-waves system.
  
By plugging \eqref{ansatz}  in \eqref{NLSd} and  by  separating real and imaginary parts,  we can rewrite \eqref{NLSd} as
 the  system
\be
\label{PhAmd}
 \left\{\begin{array}{ll}
\displaystyle{\e^2 c \p_t A^\e - c \p_x A^\e 
+ \e^2 \p_x A^\e \p_x \vp^\e 
+ \frac{1}{2} \big( 1+\e^2 A^\e \big) \p^2_x \vp^\e 
+ \e^4 \nabla_\perp A^\e \cdot \nabla_\perp \vp^\e } \\ \hspace{7.95cm}
\displaystyle{+ \frac{\e^2}{2} \big( 1+\e^2 A^\e \big) \Delta_\perp \vp^\e = 0 }
\\ \ \\ 
\ds{\e^2 c \p_t \vp^\e - c \p_x \vp^\e 
- \e^2 \frac{\p^2_x A^\e }{2\big( 1 + \e^2 A^\e \big)} 
- \e^4 \frac{\Delta_\perp A^\e }{2\big( 1 + \e^2 A^\e \big)}
+ \frac{\e^2}{2}\big( \p_x \vp^\e \big)^2 
+ \frac{\e^4}{2} |\nabla_\perp \vp^\e |^2} \\ \hspace{7.95cm}
\displaystyle{+ \frac{1}{\e^2}\, f\big( (1+\e^2 A^\e)^2 \big)} = 0.
\end{array}\right.
\ee
Now, assuming that  $A^\e \to A$ and $\vp^\e \to \vp$ as $\e \to 0$, we 
formally obtain from the two equations of the  above system that
\be
\label{deriv}
- c \p_x A + \frac{1}{2}\, \p^2_x \vp =0, \quad \quad \quad 
- c \p_x \vp + 2f'(1) A =0.
\ee
Note that   we  have used  that $f(1)=0$ and  thus that  $f\big( (1+\e^2 A^\e)^2 \big) \simeq 2\e^2 f'(1) A$
at leading order. In \iref{deriv} and from the definition 
\iref{defcoeff} of $c$, the first equation is just $- \ds{\frac{1}{2c}}$ 
times the derivative of the second equation with respect to $x$,
hence, we have found for the limit the constraint 
\be
\label{cons1} 2 c A = 
\p_x \vp.
\ee
To get the limit equation  satisfied by $A$, we can
 add   the first equation in   \iref{PhAmd} and 
$\ds{\frac{1}{2c}}$ times the derivative of the second equation with respect 
to $x$ in order to cancel the most singular term.  This yields the  equation
\be
\label{PhAm2d}
\quad \left. \begin{array}{ll}
\displaystyle{ c \p_t \big( A^\e + \frac{1}{2c}\, \p_x \vp^\e \big) 
- \frac{1}{4c} \p_x \Big( \frac{\p^2_x A^\e}{1 + \e^2 A^\e} \Big) 
+ \frac{1}{2} \big( 1 + \e^2 A^\e \big) \Delta_\perp \vp^\e 
+ \frac{c}{\e^4}\, \p_x \big( Q(\e^2 A^\e) \big) } 
\\ \ \\ \quad \quad \quad \quad 
\displaystyle{ + 
\Big\{ \p_x A^\e \p_x \vp^\e + \frac{1}{2}\, A^\e \p^2_x \vp^\e + 
\frac{1}{4c} \, \p_x \big( ( \p_x \vp^\e )^2 \big) + 
\frac{1}{2c} \big[ f'(1) + 2 f''(1) \big] \p_x \big( (A^\e)^2 \big) \Big\} } 
\\ \ \\ \quad \quad \quad \quad \quad \quad \quad \quad 
\displaystyle{ = \frac{\e^2}{4c} \p_x \Big( \frac{\Delta_\perp A^\e }{\big( 1 + \e^2 A^\e \big)} \Big) - \frac{\e^2}{4c} \p_x \big( |\nabla_\perp \vp^\e|^2 \big) }
\end{array}\right.
\ee
where
$$ c^2 Q(r) \equiv f\big( (1+r)^2 \big) - 2 f'(1) r - 
\big( f'(1) + 2 f''(1) \big) r^2 = \BO(r^3) \quad r \to 0.$$
Still on  a formal level, if  $A^\e \to A$ and $\vp^\e \to \vp$ as $\e \to 0$,  this yields 
$$ 2 \p_t A + \Big[ 6 + \frac{2}{c^2} f''(1) \Big] A \p_x A - 
\frac{1}{4c^2} \p^3_x A + \frac{1}{2c}\, \Delta_\perp \vp = 0$$
by using the relation \eqref{cons1}.  Consequently, we have obtained the sytem 
\be
\label{systKP}
\quad \quad \left\{\begin{array}{ll}
\displaystyle{\p_x \vp = 2 c A }\\ \ \\
\displaystyle{ 2 \p_t A + \Big[ 6 + \frac{2}{c^2} f''(1) \Big] A \p_x A - 
\frac{1}{4c^2} \p^3_x A + \frac{1}{2c}\, \Delta_\perp \vp = 0}
\end{array}\right.
\ee
which is a reformulation of the KP-I equation. Note that in dimension $1$, {\it i.e.} when
 $n=1$, this amounts to assume that  all the functions involved in the derivation do not
 depend on $x_{\perp}$, then the equation for $A$ in \eqref{systKP} just reduces to  
  the KdV equation since $\Delta_{\perp}\varphi=0$.

Finally, let us  notice that because of the scaling \eqref{scaling},  for the solution $\Psi$ of  the original (NLS) equation with time-scale $1$,   the 
convergence to KdV or KP-I dynamics takes place for times of order $\eps^{-3}$.\\

In dimension $n=1$, the   formal derivation of the KdV equation from the (NLS) equation in this asymptotic regime 
is well-known in the physics literature (see, for example, \cite{KAL}), and is useful 
in the stability analysis of dark solitons or travelling waves of small energy. 
In the case of  the Gross-Pitaevskii equation, for instance (that is for $f(R)=R-1$), 
the travelling waves are solutions to  (NLS)  of the form 
$\Psi(\tau,z) = U(z- \sigma \tau)$, so that $U$ solves
\be
\label{OP} -i\sigma \p_z U + \frac{1}{2} \p_{zz} U = U ( |U|^2 - 1 ), \quad \quad z\in \R\ee
with the condition $|U|(z) \to 1$ as $z \to \pm \ii$. For this nonlinearity, explicit 
integration (see, {\it e.g.} \cite{Tsu}) gives  for $0< \sigma < 1$ the nontrivial solution
$$ U_\sigma(z) = \sigma -i \sqrt{ 1 - \sigma^2 } \, {\rm th} 
\Big( z \sqrt{1 - \sigma^2} \Big) .$$
The small energy regime corresponds to $\sigma \simeq  1$, thus we set 
$\sigma^2 = 1- \e^2$, $\e>0$ small, and we obtain
$$ U_\sigma(z) = -i \e\, {\rm th} ( \e z ) + \sqrt{1 - \e^2} 
= \sqrt{1 - \frac{\e^2}{{\rm ch}^2 (\e z)}} \exp \big( i\e \vp^\e(\e z) \big), $$
with $\vp^\e(\e z) = - {\rm th} ( \e z ) + \BO(\e^3)$, and we see that this 
corresponds to the ansatz \iref{ansatz} as $\e \to 0$. Furthermore, here, 
$A^\e =  -1/\mbox{ch}^2 $ does not depend on $\e$ and is the soliton of 
the KdV equation ($c=1$, $k=6$). Note that \eqref{OP} is also often adimensionalized 
in the form
$$ - i\sigma \p_z U + \p_{zz} U = U ( |U|^2 - 1 ).$$
In this case the critical speed one, the speed of sound, is changed for $\sqrt{2}$.

In higher dimensions $n=2$, $3$, the convergence of the travelling waves to the 
Gross-Pitaevskii equation ({\it i.e.}  (NLS)  with $f(R)=R-1$) with speed 
$\simeq 1$ to a soliton of the KP-I equation is formally derived in the 
paper \cite{JR}, while in \cite{BR}, this KP-I asymptotic regime for (NLS) in 
dimension $n=3$ is used to investigate the linear  instability of the solitary waves 
of speed $\simeq 1$. On the mathematical level, in dimension $n=2$, the 
convergence of the travelling waves of speed $\simeq 1$ for the 
Gross-Pitaevskii equation to a ground state of the KP-I equation is proved 
in \cite{BGS}.

Here we shall study  the rigorous derivation  of KdV/KP-I from (NLS)
 for arbitrary time dependent solutions. All our results are in 
particular valid for the Gross-Pitaevskii equation $f(R)= R-1$.
 
In arbitrary dimension, we shall justify the KdV/KP-I  limit
   by studying directly  an hydrodynamical formulation of \eqref{NLSd} as a singular PDE limit
    as in \cite{Klainerman-Majda}, \cite{GrenierS}, \cite{Schochet}:
     we shall first prove the existence of $H^s$ solutions for  \eqref{PhAmd}
     with $s$ sufficiently large  on an interval of time independent of $\eps$
   and then pass to the limit by a weak compactness argument.
  Thanks to      the properties of the singular operator in \eqref{PhAmd}, we are able to pass
       to the limit for general initial data (i.e. "ill-prepared" data in the terminology
        of singular PDE limit), we need not assume that
       $2c A^\eps - \p_{x} \varphi^\e$ tends to zero at the initial time in 
order to be compatible with the constraint \eqref{cons1}.
        
    When $n=1$,  we will be able to pass to the limit directly from the global solution
     of (NLS) in the energy space towards the  solution of KdV in the energy space
       without assuming additional  regularity of
      the initial data but  with the assumption that the initial data are well-prepared
       in the sense that
       $ ||\partial_{x} \varphi_{0}^\eps - 2 c A_{0}^\eps ||_{L^2}/\eps$ tends to zero.

\subsection{KdV asymptotic regime for (NLS) in the energy space}

\ \indent We  first focus on the description of our result in  the one dimensional 
case $n=1$, and work only in the energy space for (NLS) and the $H^1$ energy space for 
KdV. The Cauchy problem for (NLS)  is not standard because 
of the condition at infinity $|\psi| \to 1$ (see  \cite{Ge}, \cite{Z}, \cite{Ga}) which is expected
 in order to give a meaning to the energy $\mathcal{E}(\Psi)$.  We have the following: 
\begin{theo}(\cite{Z})
\label{CauchyNLS}
There exists $\mathcal{E}_{0} >0$ such that  
for every  $\Psi_0 \in H^1_{loc}(\R)$ verifying $\BE(\Psi_0) \leq \mathcal{E}_{0}, $ 
and $ |\Psi_0|(z) \to 1$ as $|z| \to +\ii$, there exists a unique 
solution $\Psi$ to (NLS) such that
 $\Psi - \Psi_0 \in \BC \big( \R_+,H^1(\R) \big)$. 
Moreover, $\BE\big( \Psi(t) \big) = \BE(\Psi_0)$ for $t\geq 0$.
\end{theo}

This Theorem is not exactly formulated under this form in \cite{Z} (Theorem III.3.1).
Nevertheless, as we shall see in Lemma \ref{loinzero},
 if $\mathcal{E}(\Psi) \leq \mathcal{E}_{0}$ is sufficiently small, then we can write 
  $\Psi= \rho e^{i \phi}$ with
  $$ |\!|\partial_{x} \rho |\!|_{L^2} + |\!| \rho- 1 |\!|_{L^\infty} 
+ |\!|\partial_{x}\phi |\!|_{L^2} $$
  sufficiently small and hence  we can indeed use \cite{Z} (Theorem III.3.1).

\bigskip

It is also known that the Cauchy problem for the $KdV$ equation\footnote{Here, it might happen that $k=0$, in 
which case the KdV equation reduces to the so-called (linear) Airy equation $2\p_t v - \frac{1}{4c^2} \p_x^3 v =0$ and the Cauchy problem is then trivial to solve.} \cite{KPV} is well-posed in the energy space:

\begin{theo} (\cite{KPV})
\label{CauchyKdVH1} We consider the Cauchy problem for the KdV equation
$$ 2 \p_t v + k \, v \p_x v - \frac{1}{4c^2}\, \p_{xxx} v = 0, \quad \quad v_{|t=0} = v_0 .$$
 If $v_0 \in H^1(\R)$, then there exists a unique solution of the KdV equation satisfying 
$v \in \BC \big( \R_+, H^1(\R) \big)$ and $\p_x v \in L^4_{loc} \big( \R_+, L^\ii(\R) \big)$. 
\end{theo}
Note that it is possible to prove the well-posedness of KdV in spaces of much lower regularity
 than $H^1$ (see \cite{KPV2} for example) but we shall not use these results  here. 

Our first  result relates the solution of (NLS) obtained in Theorem \ref{CauchyNLS}
 in the scaling \eqref{scaling}
 and the solution of KdV obtained in Theorem \ref{CauchyKdVH1}: 

\begin{theo}
\label{KdVH1} $\bs{(n=1)}$ 
Assume  that   $(A^\e_0)_{0 < \e < 1} \in H^1$ and 
$(\vp^\e_0)_{0 < \e < 1} \in \dot{H}^1$ 
enjoy the uniform estimate
\be
\label{Mbound}
M \equiv \sup_{0 < \e < 1} \Big\{ \big|\!\big| A^\e_0 \big|\!\big|_{H^1} + 
\frac{1}{\e} \big|\!\big| \p_x \vp^\e_0 -2c A^\e_0 \big|\!\big|_{L^2} \Big\} < +\ii \ee
and that
$$ A^\e_0 \to A_0 \quad {\it in} \quad L^2 \quad {\it as} \ \ \ \e \to 0. $$
Consider the initial datum
\be
\label{psi0}
\psi^\e_0 = \big( 1 + \e^2 A^\e_0 \big) \exp\big( i\e \vp^\e_0 \big) 
\ee
for  \iref{NLSd}, and let $\psi^\e \in \psi^\e_0 + \BC ( \R_+ , H^1)$ be the associated 
solution to \iref{NLSd}(given by Theorem \ref{CauchyNLS}). 

Then, there exists $\e_0 > 0$, depending only on $M$, such that, for $0 < \e \leq \e_0 $, there 
exist two real-valued functions $\vp^\e$, $A^\e \in \BC( \R_+ \times \R, \R)$ such that 
$(A^\e,\vp^\e)_{|t=0} = ( A^\e_0,\vp^\e_0)$,
 and 
 \be
 \label{modarg}
 \psi^\e = \big( 1 + \e^2 A^\e \big) \exp\big( i\e \vp^\e \big)
 \ee
  with $1 + \e^2 A^\e \geq \frac{1}{2}. $
Furthermore, as $\e \to 0$, we have the convergence
$$  A^\e \to A\quad \mbox{ in }  \quad \mathcal{C}([0, T], H^s), \quad  \p_x \vp^\e \to 2c A, \quad \mbox{ in }  \quad \mathcal{C}([0, T], L^2)$$
for every $s <1$ and every $T>0$, where $A$ is the solution of KdV with initial value $A_{0}$.
 \end{theo}

Let us emphasize that the initial data are well-prepared (see \eqref{cons1}) in the sense that
\be
\label{as1}  \big|\!\big| \p_x \vp^\e_0 -2c A^\e_0 \big|\!\big|_{L^2} = \BO(\e). \ee
Under a stronger assumption on the preparedness of the initial data,  namely
\be
\label{as2}
 \big|\!\big| \p_x \vp^\e_0 -2c A^\e_0 \big|\!\big|_{L^2} = o(\eps),
 \ee
one can reach  the convergence in $H^1$ for the amplitude  (see Theorem  \ref{conver} in 
Subsect. \ref{H1+}). This assumption will not be needed when we work with more regular data as in 
Theorem \ref{Smoothy} below.   Finally, note that
the usual assumption of well-prepared data for a  singular system  (see \cite{Klainerman-Majda}
for example) like  \eqref{PhAmd}  in order  to get  that $\partial_{t}A^\eps= \mathcal{O}(1)$
would be that 
$$ \big|\!\big| \p_x \vp^\e_0 -2c A^\e_0 \big|\!\big|_{L^2} = \mathcal{O}(\eps^2).$$
 Consequently, we note that our assumptions 
  \eqref{as1} and even  \eqref{as2} are weaker.  
   
 Related  results are obtained in \cite{BGSS} for the Gross-Pitaevskii equation ($f(R)=R-1$)
  by using different methods, namely the complete integrability of the equation through the
  conservation of higher order energies. 

\bigskip

  The strategy of the proof is as follows. By using the conservation of the energy and
   of the momentum 
   $$ \BP = \frac{1}{2} \int_\R \big( i \Psi , \p_z \Psi \big) \ dz$$
   (actually one of its variants since $\BP$ is not well-defined for
    functions which tend to $1$ at infinity),
     we shall prove that  one can write 
     $$\psi^\e = \big( 1 + \e^2 A^\e \big) \exp\big( i\e \vp^\e \big), \quad \quad \quad \quad$$ 
 with $1 + \e^2 A^\e \geq \frac{1}{2}$ 
and the uniform bounds 
$$ \sup_{0 < \e < \e_0, \ t\in \R_+ } \Big\{ \big|\!\big| A^\e \big|\!\big|_{H^1} 
+ \frac{1}{\e} \big|\!\big| \p_x \vp^\e -2 c A^\e \big|\!\big|_{L^2} \Big\} < +\ii .$$
The $H^1$ bound  on $A^\eps$ will provide compactness in space.
  Then we shall get compactness in time by using the properties
   of the singular  part of the equation \eqref{PhAmd} namely properties
    of the transport equation with high speeds
   \be
   \label{transportintro}
   \left\{\begin{array}{l}
\ds{\p_{t} A^\eps - \frac{1}{\eps^2} \p_{x} \big( A^\e - u^\e \big) = S_A^\e} \\ \\
\ds{\p_{t} u^\eps - \frac{1}{\eps^2} \p_{x} \big( u^\e - A^\e \big) = S_u^\e}.
\end{array}\right.
\ee
This will allow to extract a subsequence which converges strongly in $L^2_{loc}(\mathbb{R}_{+}
 \times \mathbb{R}) $  towards
 the solution of the KdV equation.   Finally we shall prove that we actually have
  a better convergence  which is in particular global in space as stated in the theorem.

\subsection{KdV and KP-I asymptotic regimes for smooth initial data}

\ \indent In arbitrary dimension, we will work with $H^s$ norms and local in time smooth 
solutions in $H^s$, with $s$ sufficiently large. 

Our first result is: 
\begin{theo}
\label{Smoothy}
Let $n \geq 1$ and let $s$ such that $s> 1 + \ds{\frac{n}{2}}$. Assume that
\be
\label{bornini}
M_s \equiv \sup_{0 < \e< 1} \big|\!\big| \big( A^\e_0, \p_x \vp^\e_0, \e \nabla_\perp \vp^\e_0\big) \big|\!\big|_{H^{s+1}} <+ \infty
\ee
and consider the initial datum for \iref{NLSd}
$$ \psi^\e_0 = \big( 1 + \e^2 A^\e_0 \big) \exp\big( i\e \vp^\e_0 \big). $$
Then, there exist $T>0$ and $0 < \e_0 < 1$, depending on $M_s$, such that, 
for $0 < \e \leq \e_0 $, there exists a unique solution $\psi^\e$ to \eqref{NLSd} 
with $\psi^\e_{|t=0} = \psi^\e_0$ such that $\psi^\e - \psi^\e_0 \in \BC\big( [0,T],H^{s+1}\big)$. 
Furthermore, there exist two real-valued functions 
$A^\e \in \BC \big( [0,T] , H^{s+1} \big) $ and $\vp^\e \in \BC \big( [0,T] , \dot{H}^{s+1} \big) \cap \BC \big( [0,T] \times \R^n)$ such that $(A^\e,\vp^\e)_{|t=0} = (A^\e_0,\vp^\e_{0}) $ 
and, for $0\leq t\leq T$,
\be
\label{polaire}
 \psi^\e = \big( 1 + \e^2 A^\e \big) \exp\big( i\e \vp^\e \big), \quad \quad \quad \quad 
1 + \e^2 A^\e \geq 1/2 \ee
and
\be
\label{unif}
\sup_{0 < \e < \e_0, \ t\in [0,T]} \Big\{ \big|\!\big| A^\e \big|\!\big|_{H^{s+ 1 }(\R^n)} 
+ \big|\!\big| \big( \p_x \vp^\e, \e \nabla_\perp \vp^\e \big) \big|\!\big|_{H^{s}(\R^n)} \Big\} <+ \infty .
\ee
\end{theo}
The important result in Theorem \ref{Smoothy}  is the qualitative information that there exists
 a uniform time $T$ for which the representation \eqref{polaire} and the uniform
  bounds \eqref{unif} hold.

 To   prove  Theorem \ref{Smoothy}
   we shall  rewrite \eqref{NLSd} as a hydrodynamical equation.
   As in \cite{Grenier}, we shall use a modified Madelung transform where we allow 
the amplitude to be complex. This allows to get an hydrodynamic system with a 
much simpler structure than \eqref{PhAmd}. It is a first order hyperbolic system   
    with a singular perturbation made of a skew-symmetric zero order term and a skew-symmetric
     second  order term.
       The uniform time existence for the obtained system will then follow from uniform $H^s$ estimates 
as in the works \cite{Klainerman-Majda}, \cite{Grenier}, \cite{Schochet}.

In the recent work \cite{BDS},  the linear wave regime for the Gross-Pitaevskii equation is investigated. 
This regime occurs for  larger data  on a shorter time.
 In this regime the equivalent of Theorem \ref{Smoothy} is obtained in \cite{BDS}.
 The proof in \cite{BDS}  is different from ours since the the uniform bounds
  are obtained through the study of a different hydrodynamical system (namely the one
  obtained by the standard Madelung transform).

The next step will be the study of the  convergence towards  solutions of the KP-I equation
 of the solutions constructed in Theorem \ref{Smoothy}.
 
 Note that  for $A_{0}$ in $H^s$ with $s>1+n/2$, 
the Cauchy problem for the KP-I equation is well-posed: there exists a unique local in 
time $H^s$ solution. Note that it is actually known to be  well-posed in spaces of much lower regularity \cite{IKT}, \cite{MST}. Moreover, in dimension $n=2$, the  solutions are global in time whereas  in dimension $n=3$, the solution of KP-I may blow-up (in $H^1$) in finite time (see \cite{L}).

Our first convergence result is: 

\begin{theo}
\label{asympt1}
Under the assumptions of Theorem \ref{Smoothy}, if moreover  
there holds
\be
\label{CVinit}
(A^\eps_{0}, \partial_{x} \varphi_{0}^\eps, \eps \nabla_{\perp}\varphi_{0}^\eps) \to (A_{0},\partial_{x} \varphi_{0}, 0) \quad \quad 
{\it in} \quad L^2.
\ee
Let $A$ be the solution of the KP-I equation
$$ \partial_{x} \Big( 2  \partial_{t} A + k A \partial_{x} A - 
\frac{ 1}{ 4c^2} \partial_{x}^3 A \Big) + \Delta_{\perp} A=0 $$
with initial value $A_{/t=0}= \ds{\frac{1}{2} \big( A_{0}+ \frac{1}{ 2 c } \p_{x}\varphi_{0} \big)} 
\in H^{s+1}$. Then, we have the weak convergences, as $\e \to 0$,
$$ A^\e \rightharpoonup A \quad \quad \p_x \vp^\e \rightharpoonup 2 c A \quad \quad 
\mbox{\it{ weakly in}} \quad  L^2\big( [0, T]\times \mathbb{R}^n \big) $$
 and  the strong convergence
$$ \frac{1}{2} \Big( A^\e + \frac{1}{2c} \p_x \vp^\e \Big) \to A \quad \quad {\it in} \quad \quad 
L^2\big( [0,T], H^\s (\R^n) \big) \quad \forall \ \s < s. $$
\end{theo}

Note that the result of Theorem \ref{asympt1} holds  for smooth but ill-prepared initial data in
 the sense that they do not satisfy the constraint \eqref{cons1}.
  We shall actually get in the proof of Theorem \ref{asympt1} a stronger type of convergence.
  Namely, we get that $\partial_{x}A^\eps$ and $ (\partial_{xx}\varphi^\eps)/ 2c$
   converge strongly to $\p_x A$ in $ L^2_{loc}(0,T, H^m_{loc})$ for every $m<s$
    if $n \geq 2$ and that $A^\eps$ and $(\partial_{x}\varphi^\eps)/2c$
   converge strongly to $A$ in $ L^2_{loc}(0,T, H^{m+1}_{loc})$ if $ n=1$.

   Finally,  for slightly  well prepared  data, we are able to recover global 
strong convergence in space:
\begin{theo}\label{H1perp}
 Under the same assumptions as in Theorem \ref{Smoothy}
  and \ref{asympt1}, i.e.  \eqref{bornini} and \eqref{CVinit},  
we assume moreover that 
\begin{align}
\label{assezprep}
(n=1) \quad \quad & \big|\!\big| \p_x \vp^\e_0 - 2c A^\e_0 \big|\!\big|_{L^2} 
\to 0 \quad {\rm as} \quad \e \to 0 \nonumber, \\
(n\geq 2) \quad \quad & \big|\!\big| \p_x \vp^\e_0 - 2c A^\e_0 \big|\!\big|_{L^2} = \mathcal{O}(\eps), \quad
 |\!|  \nabla_{\perp} \varphi_{0}^\eps |\!|_{L^2} = \mathcal{O}(1).
\end{align}
Then, we have the convergences, as $\e \to 0$,
$$A^\eps \rightarrow A\quad  \mbox{ \it strongly in } \mathcal{C}([0, T ], H^m), \quad \quad \partial_{x} \varphi^\eps\rightarrow 2cA\quad  \mbox{\it  strongly  in } \mathcal{C}([0, T ], H^{m-1})$$
for every $m<s+1 $. Furthermore,  if $n\geq 2$, there exists  $K>0$  such that, 
for $0\leq t \leq T$, $0< \e < \e_0$,
\be
\label{gradperp}
\int_{\R^n} |\nabla_\perp \vp^\e|^2\ dX \leq K.
\ee
\end{theo}

We emphasize that in dimensions $n\geq 2$, the hypothesis in the last theorem  is stronger than in dimension $n=1$ 
in order to ensure the bound for $\ds{\int_{\R^n} |\nabla_\perp \vp^\e|^2\ dX }$. Moreover, 
in dimension $n=1$, \iref{assezprep} is weaker than the hypothesis in Theorem \ref{KdVH1}.

    
The paper is organized as follows.  Section \ref{section1} is devoted to the proof of 
Theorem \ref{KdVH1}, section \ref{section2} is devoted  to the Proof of Theorem \ref{Smoothy}. 
The proofs of Theorems \ref{asympt1},  \ref{H1perp} are  finally  given in sections  \ref{sectionas}, \ref{section3}.


\section{Proof of Theorem \ref{KdVH1}}
\label{section1}

\ \indent We shall split the proof in many steps.
In the first step we   prove that the modulus of  a  solution of  (NLS)  remains away from zero
 if its energy is sufficiently small
 so that it can be written as \eqref{modarg} and we prove that one can define a variant
  of the momentum which is well-defined.
   Then we shall use the energy and the momentum to get uniform $H^1\times \dot{H}^1$
    estimates for $(A^\eps, \varphi^\eps)$.
   The third step will be the study of the system \eqref{transportintro} in order to get
    compactness in time.
     Finally, the last part will be devoted to the passage to the limit in the equation.
     
     \subsection{Preliminaries}

\ \indent For the regime of 
interest  to  us, the energy is small. In this case,  we shall prove that the modulus $|\Psi|$  remains 
close to $1$.  A first useful remark is that  since $F'(1) = 2f(1) = 0$ and $F''(1) = 2f'(1)= 2c^2 >0$, 
we have for some $\delta \in (0, 1/2)$
\be
\label{below}
F(R) \geq \frac{c^2}{2} (R - 1)^2, \quad \quad | R -1 | \leq \delta
\ee
and also
\be
\label{above}
F(R) \leq C (R - 1)^2, \quad \quad | R -1 | \leq \delta
\ee
for some $C>0$.

\begin{lem}
\label{loinzero}
There exists $\BE_0 > 0$, depending only on the nonlinearity $f$, such 
that if $\Psi \in H^1_{loc}(\R)$ verifies $\BE(\Psi) < \BE_0$ and 
$|\Psi|(z) \to 1$ for $z \to +\ii$, then
$$ \big|\!\big| \, |\Psi|^2 -1 \big|\!\big|_{L^\ii(\R)} \leq \delta.$$
\end{lem}

Note that for  an initial value under the form
 \eqref{psi0}, we have,  since $M$ is finite,  that 
 $$ \mathcal{E}(\psi_{0}^\eps) \leq C \eps \Bigl( \int_{\mathbb{R}}  \eps^4 ( \partial_{x} A_{0}^\eps )^2
  +  \eps^2( 1 + \eps^2 A_{0}^\eps)^2 (\partial_{x} \varphi_{0}^\eps )^2
   + \eps^2 (A_{0}^\eps)^2   \, dx\Big) \leq C \eps^3$$
   where $C$ depends only on $M$. Consequently,  since the energy is conserved, 
    we can indeed use Lemma 
 \ref{loinzero}
    for $\eps$ sufficiently small  to write the solution  $\psi^\eps$ of NLS
     given by Theorem \ref{CauchyNLS} under the form  $\psi^\eps= \rho e^{i \phi}$
     with $\phi\in H^1_{loc}$  and $|\rho^2 -1 |\leq 1/2$. Note that $\rho$ and $\phi$ depend on $\eps$ but we
      omit this dependence in our notation.\\

\noindent {\bf Proof of Lemma \ref{loinzero}.} Since $|\Psi|(z) \to 1$ for $z \to +\ii$, 
we have
\be
\label{primit}
\forall z\in \R, \quad \quad \big( |\Psi|^2(z) - 1 \big)^2 = 
- 4 \int_z^{+\ii} |\Psi| \big( |\Psi|^2 - 1 \big) \p_z |\Psi|,
\ee
and we can define the maximal interval $I=[a,+\ii)$ such that 
$\big| |\Psi|^2 -1 \big| \leq \delta$ in $I$. Then,
$$ \int_I \big| \p_z |\Psi| \big| ^2 + \frac{c^2}{2} \big( |\Psi| - 1 \big)^2 \ dz \leq 
\int_\R |\p_z \Psi|^2 + F \big( |\Psi|^2 \big) \ dz = 
2 \BE(\Psi).$$
As a consequence, by \iref{primit} and Cauchy-Schwarz,
$$ \big|\!\big| \, |\Psi|^2 -1 \big|\!\big|_{L^\ii(I)}^2 \leq 
4 \sqrt{ 1 + \delta}  |\!| \, |\Psi^2| -1 |\!|_{L^2(I)} \cdot 
\big|\!\big| \p_z |\Psi| \big|\!\big|_{L^2(I)} \leq K_0 \BE(\Psi) ,$$
where $K_0$ depends only on $f$.  The result follows from an easy 
continuation argument, taking $\BE_0 \equiv \delta^2/ K_0$.\b \\

 Next, we recall that the Schr\"odinger flow  also formally preserves the {\it momentum}, 
that should be defined by
$$ \BP = \frac{1}{2} \int_\R \big( i \Psi , \p_z \Psi \big) \ dz.$$
However, this quantity does not make sense as a Lebesgue integral for a map $\Psi$ which is 
just of finite energy with $|\Psi| \to 1$ at infinity. Notice that if 
$\Psi= \rho \exp(i\phi)$, then
$$ \BP =  \frac{1}{2} \int_\R \rho^2 \p_z \phi \ dz.$$
Variants of the momentum $\BP$ are also formally conserved by the 
Schr\"odinger equation (NLS), 
namely
$$ \frac{1}{2} \int_\R \big( i (\Psi - 1) , \p_z \Psi \big) \ dz \quad \quad {\rm if} \ \ 
\Psi \to 1 \quad {\rm at \ infinity}$$
and
$$ \frac{1}{2} \int_\R \big( \rho^2 -1 \big) \p_z \phi \ dz.$$
This last integral has the advantage to be a Lebesgue integral if $ \Psi \in H^1_{loc}(\R)$ satisfies 
$$ \BE(\Psi) < +\ii, \ \ |\Psi|(x) \to 1 \ \ {\rm as} \ x \to +\ii\ \quad {\rm and} \ \quad \big|\,
|\Psi|^2 - 1 \big| \leq \delta, $$
since then
$$ \frac{1}{2} \int_\R \big( \rho^2 -1 \big) \p_z \phi \ dz = 
\frac{1}{2} \int_\R \big( |\Psi|^2 -1 \big) \frac{ {\rm Im}( \p_z \Psi)}{|\Psi|} \ dz.$$
As we have seen in the remark  after  Lemma \ref{loinzero}, in our regime, the map $\psi^\eps $  satisfy the bound 
$\big|\!\big| \, |\psi^\eps|^2 -1 \big|\!\big|_{L^\ii(\R)} \leq \delta $ and  hence, we  have a 
well-defined momentum, if we take this last definition.

 Finally, in view of the scaling \eqref{scaling}, it is usefull to introduce a rescaled version of
  the energy. We set
\begin{eqnarray}
\nonumber
 E^\eps(\psi^\eps) &=& \frac{\mathcal{E}(\Psi)}{ \eps } = \frac{ 1 }{2} \int_{\R} |\p_{x} \psi^\eps |^2
 + \frac{1}{\eps^2} F( |\psi^\eps|^2) \, dx\\
 \label{energie} & = & \frac{ 1}{2} \int_{\R} |\p_{x} \rho|^2 + \rho^2 |\p_{x}\phi|^2 + 
\frac{1}{ \eps^2}F(\rho^2) \,dx
 \end{eqnarray}
 since  $\psi^\eps= \rho e^{i\varphi}.$ In a similar way, we define a rescaled momentum
\be
\label{rescmom} P^\e(\psi^\eps ) \equiv \frac{\e}{2} \int_\R \big( \rho^2 -1 \big) \p_x \phi \ dx.
\ee 
Note that both quantities are conserved.

\subsection{Uniform estimates}

\ \indent We shall prove the following:

\begin{lem}\label{Bornes} Under the assumptions of Theorem \ref{KdVH1}, there 
exists $\e_0 > 0$, depending only on $M$, such that, for $0 < \e \leq \e_0 $, there 
exist two real-valued functions $\vp^\e$, $A^\e \in \BC( \R_+ \times \R, \R)$ such that 
$(A^\e,\vp^\e)_{|t=0} = ( A^\e_0,\vp^\e_0)$,
$$\psi^\e = \big( 1 + \e^2 A^\e \big) \exp\big( i\e \vp^\e \big), \quad \quad \quad \quad 
1 + \e^2 A^\e \geq \frac{1}{2} $$
and
$$ \sup_{0 < \e < \e_0, \ t\in \R_+ } \Big\{ \big|\!\big| A^\e \big|\!\big|_{H^1} 
+ \frac{1}{\e} \big|\!\big| \p_x \vp^\e -2 c A^\e \big|\!\big|_{L^2} \Big\} < +\ii .$$
\end{lem}

\subsection*{Proof of Lemma \ref{Bornes}.}

\ \indent 
The proof relies on the use of the  conservation of $E^\eps$ and $P^\eps$ as noticed 
in \cite{BGS}. In particular, the quantity $E^\eps -2c P^\eps$ gives valuable information.

As we have already seen, we can write    $\psi^\eps = \rho \exp(i\phi)$ for 
some real-valued functions $\rho \geq 1/2$ and $\phi$ in $H^1_{loc}(\R)$.
Note that 
$$ |\p_x \psi^\e|^2 = (\p_x \rho)^2 + \rho^2 (\p_x \phi)^2.$$ 
Next, we set
$$ F (R) = c^2 \big( R -1 \big)^2 + F_3(R ), \quad {\rm with} \quad 
F_3(1 + r ) = \BO(r^3), \quad r\to 0.$$
 By using \eqref{energie} and \eqref{rescmom}, this yields  
\be
\label{develo}
E^\e(\psi^\e) = \frac{1}{2} \int_\R (\p_x \phi)^2 
+ \frac{c^2}{\e^2} \big( \rho^2 - 1 \big)^2 
+ (\rho^2 - 1) \cdot (\p_x \phi)^2 + (\p_x \rho)^2 
+ \frac{1}{\e^2}\, F_3 ( \rho^2 - 1) \ dx
\ee
and
\be
\label{develoE-cP}
E^\e(\psi^\e) - 2c P^\e(\psi^\e) = \frac{1}{2} \int_\R 
\big( \rho^2 -1 \big) (\p_x \phi)^2 + (\p_x \rho)^2 + 
\Big( \p_x \phi - \frac{c}{\e} (\rho^2-1 ) \Big)^2
+ \frac{1}{\e^2}\, F_3 \big( \rho^2-1 \big) \ dx,
\ee
where we have used the identity
$$ (\p_x \phi)^2 + \frac{c^2}{\e^2}\, (\rho^2-1)^2 - \frac{2c}{\e} (\rho^2-1 ) \p_x \phi = 
\Big( \p_x \phi - \frac{c}{\e} (\rho^2-1 ) \Big)^2. $$

\bigskip

The proof  of Lemma \ref{Bornes} is divided in 3 Steps. In the proof, $K$ stands for a constant 
depending only on $f$ and $M$. \\

\noindent {\bf Step 1:} We first prove  the following expansions for $E^\e(\psi_0^\e)$ 
and $E^\e(\psi_0^\e) - 2c P^\e(\psi_0^\e)$ as $\e \to 0$:
$$ E^\e (\psi^\e_0) = 
\frac{\e^2}{2} \int_\R 4 c^2 (A_0^\e)^2 + (\p_x \vp_0^\e)^2 \ dx + \BO( \e^4) 
= 4 c^2 \e^2 \int_\R A_0^2 \ dx + o(\e^2) +  \BO( \e^4)$$
and
$$ E^\e(\psi^\e_0) - 2c P^\e(\psi^\e_0) \leq K \e^4.$$

\bigskip

This follows from  \iref{develo} and \iref{develoE-cP} with $\rho = 1+ \e^2 A^\e_0$ and $\phi = \e \vp_0^\e$. 
Indeed, from the uniform bound in $H^1$ for $A^\e_0$, we immediately infer by Sobolev 
embedding $H^1(\R) \subset L^\ii(\R)$ that 
$|\!| A^\e_0 |\!|_{L^\ii} \leq K$ and $|\psi^\e_0| = |1+\e^2 A^\e_0| \in [ 1/2, 2]$ 
for $0 < \e < \e_0$ sufficiently small, depending on $M$. Moreover, 
$\rho^2 -1 = 2 \e^2 A^\e_0 + \BO_{L^\ii(\R)}(\e^4)$. Since 
$ |F_3(R)| \leq K |R-1|^3$ for $0\leq R \leq 2$, we have 
$| F_3(\rho^2 -1) | \leq K \e^6 ( A^\e_0 )^2$, and the expansion for the energy follows. 
Concerning the expansion for $E_\e(\psi^\e_0) - 2c P_\e(\psi^\e_0)$, it suffices to use the 
assumption $|\!| \p_x \vp^\e_0 - 2 c A^\e_0 |\!|_{L^2}^2 \leq M^2 \e^2 $.\\

\bigskip

\noindent {\bf Step 2:}  We shall prove that  for every $t \in \R_+$, 
$$ \big|\!\big| \rho^2 -1 \big|\!\big|_{L^\ii(\R)} \leq K \e^2.$$

\bigskip

This will be a consequence of  the conservation of energy and momentum. Let $t\in \R_+$. We 
first infer from \iref{develoE-cP} a better estimate for 
$\ds{\int_\R (\p_x \rho)^2 \ dx}$. Since $\rho \geq 1/2$, we have, on the one hand,
\be
\label{modu1}
\Big| \int_\R (\rho^2-1) (\p_x \phi)^2 \ dx \Big| \leq 
4 \big|\!\big| \rho^2 - 1 \big|\!\big|_{L^\ii(\R)} \int_\R \rho^2(\p_x \phi)^2 \ dx 
\leq K \e^2 \big|\!\big| \rho^2 - 1 \big|\!\big|_{L^\ii(\R)},
\ee
and on the other hand, in view of $| \rho^2 -1|\leq \delta$, 
$F_3(r)=\BO(r^3)$ as $r \to 0$ there holds
\be
\label{modu2}
\Big| \int_\R \frac{1}{\e^2}\, F_3\big( \rho^2-1 \big)  \ dx \Big| 
\leq K \big|\!\big| \rho^2 - 1 \big|\!\big|_{L^\ii(\R)} 
         \int_\R \frac{1}{\e^2}\, \big( \rho^2-1 \big)^2 \ dx 
\leq K \e^2 \big|\!\big| \rho^2 - 1 \big|\!\big|_{L^\ii(\R)}.
\ee
Since $E^\e$ and $P^\e$ do not depend on time, inserting \iref{modu1} and 
\iref{modu2} into \iref{develoE-cP} yields
\begin{align*}
K \e^4 \geq E^\e(\psi^\e) - 2c P^\e(\psi^\e) \geq & \ 
\frac{1}{2} \int_\R (\p_x \rho)^2 \ dx - 
\Big| \int_\R (\rho^2-1) (\p_x \phi)^2 \ dx \Big| - 
\Big| \int_\R \frac{1}{\e^2}\, F_3\big( \rho^2-1 \big) \ dx \Big|\\
\geq & \ \frac{1}{2} \int_\R (\p_x \rho)^2 \ dx - 
K \e^2 \big|\!\big| \rho^2 - 1 \big|\!\big|_{L^\ii(\R)},
\end{align*}
so that
\be
\label{modu}
\int_\R (\p_x \rho)^2 \ dx \leq K \e^4 + 
K \e^2 \big|\!\big| \rho^2 - 1 \big|\!\big|_{L^\ii(\R)}.
\ee

We now write, since $\rho = |\psi^\e| \to 1$ as $|x| \to +\ii$,
$$ \big( \rho^2 - 1 \big)^2 (x) = 
- 4 \int_x^{+\ii} \rho \big( \rho^2 - 1 \big) \p_x \rho 
\leq C \e \sqrt{E_\e(\psi^\e)} \Big( \int_\R (\p_x \rho)^2 \ dx \Big)^{1/2} $$
by Cauchy-Schwarz inequality. From the above estimate \iref{modu} and letting
$$ \eta_\e \equiv \frac{1}{\e^2} 
\big|\!\big| \rho^2 -1 \big|\!\big|_{L^\ii(\R)},$$
we obtain
$$ \e^4 \eta_\e^2 \leq K \e^2 \sqrt{\e^4 + \e^4 \eta_\e},$$
that is
$$ \eta_\e^2 \leq K \sqrt{1 + \eta_\e}.$$
This estimate provides immediately the result
$$ \eta_\e = \frac{1}{\e^2} 
\big|\!\big| \rho^2 -1 \big|\!\big|_{L^\ii(\R)} \leq K.$$

\bigskip
\bigskip

We then set 
$$ A^\e \equiv \frac{1}{\e^2} \, ( \rho -1) 
\quad \quad {\rm and} \quad \quad 
\vp^\e \equiv \frac{\phi}{\e}. $$

\bigskip

\noindent {\bf Step 3:}  We finally prove that 
\be
\label{bonnar}
|\!| A^\e |\!|_{H^1(\R)} \leq K, \quad \quad
|\!| \p_x \vp^\e |\!|_{L^2(\R)} \leq K 
\quad \quad {\rm and} \quad \quad
|\!| 2c A^\e - \p_x \vp^\e |\!|_{L^2(\R)} \leq K \e.
\ee

\bigskip

Indeed, from Step 2, \iref{modu1} and \iref{modu2} imply
$$ \Big| \int_\R (\rho^2-1) (\p_x \phi)^2 \ dx \Big| \leq K\e^4 
\quad \quad \quad {\rm and} \quad \quad \quad 
\Big| \int_\R \frac{1}{\e^2}\, F_3\big( \rho^2-1 \big) \ dx \Big| \leq K\e^4.$$ 
Inserting this into \iref{develoE-cP} gives
$$ \Big|\!\Big| \frac{1}{\e^2} ( \rho^2 -1) \Big|\!\Big|_{H^1(\R)} \leq K, \quad \quad 
\int_\R (\p_x \phi)^2 \ dx \leq K\e^2 \quad \quad {\rm and} \quad \quad 
\int_\R \Big( \p_x \phi - \frac{c}{\e} (\rho^2-1 ) \Big)^2 \ dx \leq K\e^4 $$
and the conclusion follows. This finishes the proof of the Lemma. \b

\subsection{Properties of the wave operator}

\ \indent
In the previous subsection, we have obtained uniform bounds which  will provide
 (local) compactness in space.  We shall try now to obtain some compactness in time.

\begin{lem}
\label{source}
Consider $(A^\eps(t,x), u^\eps(t,x))$ a solution of the system
\be
\label{wave}
\left\{\begin{array}{l}
\ds{\p_{t} A^\eps - \frac{1}{\eps^2} \p_{x} \big( A^\e - u^\e \big) = S_A^\e} \\ \\
\ds{\p_{t} u^\eps - \frac{1}{\eps^2} \p_{x} \big( u^\e - A^\e \big) = S_u^\e},
\end{array}\right.
\ee
with initial data
$$ A^\e_{|t=0} = A_{0}^\e, \quad u^\e_{|t=0} = u_{0}^\e $$
and assume that, for some $\sigma \in \N$,
\begin{itemize}
\item[i)] $(A_{0}^\eps)_{0<\e<1}$ and $(u_{0}^\eps)_{0<\e<1}$ are uniformly bounded in $L^2$;
\item[ii)] $(S_A^\eps)_{0<\e<1}$ and $(S_u^\eps)_{0<\e<1}$ are uniformly bounded in 
$L^\infty\big( \mathbb{R}_{+}, H^{-\s} (\R) \big)$.
\end{itemize}
Then, for every $T>0$, $R>0$,\\

\centerline{$ (A^\eps)_{0<\e<1}$ and $ (u^\eps)_{0<\e<1}$ are uniformly bounded in 
$H^{\frac{1}{2}} \big( [0, T], H^{-\s-1}(-R,R)\big)$.}
\end{lem}

\subsection*{Proof of Lemma \ref{source}.}

These bounds come from the fact that the speed $\ds{\frac{1}{\e^2}}$ of the characteristics 
of the transport equation is extremely large compared to the size of the space domain $(-R,R)$.

 We start the proof of Lemma \ref{source} with the following lemma, where 
we take into account only the initial data, and not the source terms.

\begin{lem}
\label{init}
Consider  $(A^\eps(t,x), u^\eps(t,x))$ a solution of the system
$$\left\{\begin{array}{l}
\ds{\p_{t} A^\eps - \frac{1}{\eps^2} \p_{x} \Big( A^\eps - u^\eps \Big) = 0} \\ \\
\ds{\p_{t} u^\eps - \frac{1}{\eps^2} \p_{x} \Big(u^\eps - A^\eps\Big) = 0},
\end{array}\right. $$
with initial data
$$ A^\eps_{|t=0} = A_{0}^\eps, \quad u^\eps_{|t=0} = u_{0}^\eps.$$
Assume that $ (A^\eps_{0})_{0<\e<1}$, $(u_{0}^\eps)_{0<\e<1}$ are uniformly bounded in 
$L^2(\mathbb{R})$. Then for every $T>0,$  $ R>0$, $A^\eps$ and $ u^\eps$ are uniformly bounded in 
$H^{\frac{1}{2}} \big( [0, T], H^{-1}(-R,R) \big)$.
\end{lem}

\noindent {\bf Proof of Lemma \ref{init}.}  At first, we notice that
$$ \partial_{t} \big(A^\eps + u^{\eps} \big) =0 $$ 
and that
\be
\label{diff}
\p_{t} \big(A^\eps - u^\eps \big) - \frac{2}{\eps^2} \p_{x} \big( A^\eps - u^\eps \big) = 0.
\ee
The resolution of these transport equations gives
$$  A^{\eps}(t,x)+u^\eps(t,x)= A^\eps_{0}(x) + u_{0}^\eps(x) $$
and
\be
\label{A-u1}
A^{\eps}(t,x) - u^\eps(t,x) = A_{0}^\eps(x+  2 \eps^{-2} t) - u_{0}^\eps(x + 2 \eps^{-2}t ).
\ee
This immediately yields that
\be
\label{A+u1}
A^\eps +u^\eps \mbox{ is uniformly bounded in } H^1 \big( 0, T, L^2( \R) \big)
\ee
and hence by continuous injection, it is in particular bounded in 
$H^{\frac12} \big( 0,T, H^{-1}(-R,R)\big)$.

Next, we shall study $A^\eps - u^\eps$. From the explicit expression \eqref{A-u1}, 
we first get that
$$ \int_{0}^T \int_{-R}^R \big| A^\eps - u^\eps \big|^2(t,x) \, dxdt
= \int_{0}^T \int_{-R}^R \big| A_{0}^\eps- u_{0}^\eps \big|^2(x+2 \eps^{-2} t) \, dxdt.$$
Consequently, by using Fubini Theorem and then changing the variable $t$ into 
$\tau =x+ 2 \eps ^{-2} t$, we get
\be
\label{Ltwo}
\int_{0}^T \int_{-R}^R \big| A^\eps - u^\eps \big|^2(t,x) \, dxdt 
\leq \frac{\eps^2}{2} \int_{-R}^R |\!| A_{0}^\eps - u_{0}^\eps |\!|_{L^2(\R)}^2 \,dx \leq C R \eps^2.
\ee
In the proof, $C$ denotes a constant depending on $R$ and the uniform bounds for 
$(A^\eps_{0})_{0<\e<1}$ and $(u_{0}^\eps)_{0<\e<1}$ in $L^2$. We have thus in particular 
proven the uniform bound
\be
\label{L21}
\big|\!\big| A^\eps - u^\eps \big|\!\big|_{L^2( 0, T, H^{-1}(-R,R))} \leq C \e.
\ee
To estimate the time derivative, it suffices to remark that \eqref{diff} yields
$$ \big|\!\big| \p_{t} \big( A^\eps- u^\eps\big)(t, \cdot) \big|\!\big|_{H^{-1}(-R,R)} 
= \frac{2}{\e^2} \big|\!\big| \p_x \big( A^\eps- u^\eps\big)(t, \cdot) \big|\!\big|_{H^{-1}(-R,R)} 
\leq \frac{2}{\e^2} \big|\!\big| \big( A^\e- u^\e \big) (t, \cdot) \big|\!\big|_{L^2(-R,R)}. $$
Hence, taking the $L^2$ norm in time and using \iref{Ltwo} gives
\be
\label{L22}
\big|\!\big| \p_t ( A^\eps- u^\eps) \big|\!\big|_{L^2((0,T), H^{-1}(-R,R))} 
\leq \frac{C}{\e} .
\ee
Interpolating in time between \iref{L21} and \iref{L22}, we deduce
\be
\label{H120}
\big|\!\big| A^\e - u^\e \big|\!\big|_{H^{\frac12 }((0,T), H^{-1}(-R,R))} \leq C.
\ee
The combination of \eqref{A+u1} and \eqref{H120} ends the proof.\b

\bigskip

We shall now give the proof of Lemma \ref{source}.
Since the system \eqref{wave} is linear, we can write its solution
as the sum of the  solution of the homogeneous system and the solution of the 
nonhomogeneous system with zero initial data. Thanks to Lemma \ref{init}, we already 
know that the first term is uniformly bounded in $H^{\frac12}(0,T, H^{-1}_{loc})$ and hence 
in $H^{\frac12}(0,T, H^{-\s-1}_{loc})$. Consequently, we can focus on the second term. This 
means that we consider the solution of \eqref{wave} with zero initial value.

We notice that
$$ \partial_{t}(A^\eps + u^\eps) = S_A^\e + S_u^\e,$$
and we recall that the initial values are zero. Hence, 
$$ \big( A^\e + u^\e \big) (t) = \int_0^t \big( S_A^\e + S_u^\e \big)(s) \ ds,$$
thus we immediately get that
\be
\label{res1}
A^\eps + u^\eps \mbox{  is uniformly bounded in } H^1 \big( 0, T, H^{-\s}(\R) \big).
\ee
Similarly, since $A^\e - u^\e$ solves
\be
\label{A-u}
\p_t \big( A^\e - u^\e \big) - \frac{2}{\e^2} \p_x \big( A^\e - u^\e \big) = S^\e_A - S^\e_u
\ee
with zero initial value, we infer
$$(A^\eps - u^\eps)(t,x)= \int_{0}^t \big( S_A^\e - S_u^\e \big)\big(s, x+ 2\eps^{-2}(t-s) \big) ds.$$
By assumption $ii)$, $S_A^\e - S_u^\e$ is uniformly bounded in $L^\ii( \R_+, H^{-\s})$, 
hence, using a standard characterization of $H^{-\s}$, $\s \in \N$, there exists 
$g^\e=(g^\e_0, g^\e_1 , ... , g^\e_{\s}) \in L^\ii \big( \R_+ , L^2(\R,\R^{\s +1}) \big)$ such that
$$ S_A^\e - S_u^\e = \sum_{j=0}^{\s} \p_x^j g^\e_j. $$
Furthermore, for any interval  $I$,
$$ \big|\!\big| \big( S_A^\e - S_u^\e \big)(t,\cdot) \big|\!\big|_{H^{-\s}(I)} 
\leq \big|\!\big| g^\e \big|\!\big|_{L^2(I,\R^{\s+1}) } \leq C. $$
Here, $C$ stands for a constant depending on $R$, $T$ and the uniform bounds for 
$(A^\eps_{0} , u_{0}^\eps )_{0<\e<1}$ and $(S_A^\e,S_u^\e)_{0<\e<1}$ in $H^{-\s}$. As a 
consequence, we get from \eqref{A-u} that
\begin{eqnarray*}
\big|\!\big| \big( A^\e - u^\e \big) (t, \cdot) \big|\!\big|_{H^{-\s}(-R,R)}^2 
& \leq & 
\Big( \int_0^t \big|\!\big| g^\e \big( s,\cdot +2\e^{-2}(t-s) \big) \big|\!\big|_{L^2(-R, R,\R^{\s +1})}\, ds \Big)^2 \\ 
& \leq & 
t \int_0^t \big|\!\big| g^\e \big( s,\cdot + 2 \eps^{-2}(t-s) \big) \big|\!\big|_{L^2(-R, R,\R^{\s +1})}^2\, ds
\end{eqnarray*}
and hence that
\begin{eqnarray*}
\int_{0}^T \big|\!\big| \big( A^\eps - u^\eps \big) (t, \cdot) \big|\!\big|_{H^{-\s}(-R,R)}^2 \, dt 
& \leq & 
T \int_0^T \int_{0}^t \int_{-R}^R \big| g^\e \big|^2 \big( s, x + 2 \eps^{-2}(t-s) \big) \, dx ds dt,
\end{eqnarray*}
which we can rewrite, by using Fubini Theorem, as:
$$ \int_{0}^T \big|\!\big| \big( A^\e - u^\e \big) (t, \cdot) \big|\!\big|_{H^{-\s}(-R,R)}^2 \, dt 
\leq T \int_{-R}^R \int_{0}^T \int_{s}^T \big| g^\e \big|^2 \big( s, x + 2 \eps^{-2}(t-s) \big) \, dt ds dx.$$
By changing $t$ into $ \tau= x+ 2\eps^{-2}(t-s)$, this yields
\begin{eqnarray*}
\int_{0}^T \big|\!\big| \big( A^\e - u^\e \big) (t, \cdot) \big|\!\big|_{H^{-\s}(-R,R)}^2 \, dt
& \leq & \frac{1}{2} T \e^2 \int_{-R}^R \int_{0}^T \big|\!\big| g^\e (s,\cdot) \big|\!\big|_{L^2
 (\mathbb{R})}^2\, ds dx \leq C \e^2.
\end{eqnarray*}
We have thus proven that
\be
\label{L2}
\big|\!\big| A^\e - u^\e \big|\!\big|_{L^2( 0, T, H^{-\s}(-R,R) )} \leq C \e,
\ee
which implies in particular that 
\be
\label{Ldeux}
\big|\!\big| A^\e - u^\e \big|\!\big|_{L^2( 0, T, H^{-\s-1}(-R,R) )} \leq C \e.
\ee
To estimate $\p_{t}(A^\e - u^\e)$, we infer from \eqref{A-u}
\begin{align*}
\big|\!\big| \p_t \big( A^\e - u^\e \big) \big|\!\big|_{H^{-\s-1}(-R,R)} 
& \ \leq \frac{2}{\e^2} \big|\!\big| \p_x \big( A^\e - u^\e \big)\big|\!\big|_{H^{-\s-1}(-R,R)} 
+ \big|\!\big| S^\e_A - S^\e_u \big|\!\big|_{H^{-\s-1}(-R,R)} \\ 
& \ \leq \frac{2}{\e^2} \big|\!\big| A^\e - u^\e \big|\!\big|_{H^{-\s}(-R,R)} + C,
\end{align*}
which yields, for $ 0 < \e < 1$ and in view of \eqref{L2},
\be
\label{ache1}
\big|\!\big| \p_t \big( A^\e - u^\e \big) \big|\!\big|_{L^2(0,T, H^{-\s-1}(-R,R))} \leq \frac{C}{\e}.
\ee
Interpolation in time between \eqref{Ldeux} and 
\eqref{ache1} yields
\beq
\label{H1/2}
\big|\!\big| A^\e - u^\e \big|\!\big|_{H^{\frac12}(0, T, H^{-\s-1}(-R,R))} \leq C. 
\eeq 
To end the proof, it suffices to combine \eqref{res1} and \eqref{H1/2}.\b

\subsection{End of the proof of Theorem \ref{KdVH1}}
\label{finpreuve}

\ \indent  Since $\rho^\e = 1 + \e^2 A^\e \geq 1/2$ 
in $\R_+ \times \R$ for $0 < \e <\eps_0$, we may then rewrite \iref{NLSd} under the form \iref{PhAmd}.
In dimension $1$,  this reads
\be
\label{PhAm}
 \quad \quad  \left\{\begin{array}{ll}
\displaystyle{\e^2 c \p_t A^\e - c \p_x A^\e + 
\e^2 \p_x A^\e \p_x \vp^\e + \frac{1}{2} \big( 1+\e^2 A^\e \big) \p_{xx} \vp^\e } = 0 
\\ \ \\ 
\displaystyle{\e^2 c \p_t \vp^\e - c \p_x \vp^\e - 
\e^2 \frac{\p_{xx}A^\e }{2\big( 1 + \e^2 A^\e \big)} + 
\frac{\e^2}{2}\big( \p_x \vp^\e \big)^2 + 
\frac{1}{\e^2}\, f\big( (1+\e^2 A^\e)^2 \big)} = 0,
\end{array}\right.
\ee
and we wish to pass to the limit as $\e \to 0$. Let us define
$$ u^\e \equiv \frac{1}{2c} \p_x \vp^\e. $$

We shall first prove  that the  functions $(A^\e)_{0 < \e < \e_0}$ and $(u^\e)_{0 < \e < \e_0}$ 
are strongly precompact in $L^2_{loc}(\R_+ \times \R)$.
Indeed, we may rewrite \eqref{PhAm} as
$$ \left\{\begin{array}{l}
\ds{\p_{t} A^\eps - \frac{1}{\eps^2} \p_{x} \big( A^\e - u^\e \big) = S_A^\e} \\ \\
\ds{\p_{t} u^\eps - \frac{1}{\eps^2} \p_{x} \big( u^\e - A^\e \big) = S_u^\e},
\end{array}\right. $$
where
$$ \left\{\begin{array}{ll}
\ds{S_A^\e \equiv } & \ds{ - 2 u^\e \p_{x} A^\e - A^\e \p_{x} u^\e } \\ \\
\ds{S_u^\e \equiv } & \ds{ - \p_{x} \big( (u^\e)^2 \big) 
+ \p_{x} \Big( \frac{\p_{xx} A^\e}{4 c^2 ( 1 + \e^2 A^\eps )} \Big) 
- \frac{1}{\e^4} \p_{x}\big( \tilde{f}(\e^2 A^\e ) \big)}
\end{array}\right. $$
and
$$ \t{f} (r) \equiv \frac{1}{c^2}f\big( (1+r)^2 \big) - 2 r = \BO(r^2) \quad \quad {\rm as} \ \ r \to 0. $$

In order to use Lemma \ref{source}, we shall prove that for some constant
 $K$ depending only on $M$, we have
\be
\label{Hmoins2}
\big|\!\big| S_A^\e \big|\!\big|_{L^\ii(H^{-2})} + \big|\!\big| S_u^\e \big|\!\big|_{L^\ii(H^{-2})} \leq K.
\ee
We first note that, if $t\in \R_+$ and $ \zeta \in \BC^\ii_c (\R)$,
\begin{align*}
\langle S_A^\e(t) , \zeta \rangle = & \ - \langle u^\e(t) \p_{x} A^\e(t) , \zeta \rangle 
+ \langle u^\e (t)A^\e(t) , \p_{x} \zeta \rangle \\ 
\leq & \ \big|\!\big| u^\e(t) \big|\!\big|_{L^2}\, \big|\!\big| \p_{x} A^\e(t) \big|\!\big|_{L^2} 
\big|\!\big| \zeta \big|\!\big|_{L^\infty} + 
\big|\!\big| u^\e (t) \big|\!\big| _{L^2} \, \big|\!\big| A^\e(t) \big|\!\big|_{L^\infty }\, 
\big|\!\big| \p_x \zeta \big|\!\big|_{L^2}.
\end{align*}
Hence, by using the embedding $H^1(\R) \subset L^\infty(\R)$ and Lemma \ref{Bornes}, we get:
$$ \big|\!\big| S_A^\e (t) \big|\!\big|_{H^{-1}(\mathbb{R})} \lesssim 
\big|\!\big| u^\e(t) \big|\!\big|_{L^2}\, \big|\!\big| A^\e(t) \big|\!\big|_{H^1} \leq K. $$
In a similar way, we have, for $t\in \R_+$ and $ \zeta \in \BC^\ii_c (\R)$,
\begin{align*}
\langle S_u^\e(t) , \zeta \rangle = & \ 
\int_{\mathbb{R}} \Big[ (u^\e)^2 + \frac{1}{\e^4} g( \e^2 A^\e) 
- \frac{\e^2 (\p_{x} A^\e )^2}{4 c^2 (1 + \e^2 A^\e)^2} \Big] \p_{x} \zeta
     +\frac{\p_{x}A^\e}{4c^2 ( 1 + \e^2 A^\e)} \partial_{xx} \zeta.\\ \leq & \  K\Big(
\Big[ \big|\!\big| u^\e (t) \big|\!\big|_{L^2} + \big|\!\big| A^\e (t) \big|\!\big|^2_{L^2} 
+ \e^2 \big|\!\big| \p_x A^\e (t) \big|\!\big|^2_{L^2} \Big] \big|\!\big| \p_x \zeta \big|\!\big|_{L^\ii} 
+ \big|\!\big| \p_x A^\e (t) \big|\!\big|_{L^2} \ \big|\!\big| \p_{xx} \zeta \big|\!\big|_{L^2}\Big),
\end{align*}
where we have used that $\t{f}(r)=\BO(r^2)$ as $r \to 0$, and $\e^2 |\!| A^\e|\!|_{L^\ii} \leq 1/2$. 
Using again the embedding $H^1(\R) \subset L^\infty(\R)$ and Lemma \ref{Bornes}, this yields, for 
$0 < \e < \e_0$,
$$ \big|\!\big| S_{u}^\eps (t) \big|\!\big|_{H^{-2}(\R)} \lesssim
||u^\eps ||_{L^2}^2 + ||A^\eps ||_{H^1}^2 + ||A^\eps ||_{H^1} \leq K. $$

Consequently, thanks to \eqref{Hmoins2} and the fact that by our assumptions, $A^\e_0$ and 
$u^\e_0$ are uniformly bounded in $L^2$, we may apply Lemma \ref{source} with $\s=2$ and deduce 
that $(A^\e)_{0 < \e <\e_0} $ and $(u^\e)_{0 < \e <\e_0}$ are uniformly bounded in 
$H^{ \frac12}_{loc}( \R_+, H^{-3}_{loc})$.  In particular, since $(A^\e)_{0 < \e <\e_0}$ is 
uniformly bounded in $L^\infty(\mathbb{R}_{+}, H^1)$ and in $H^{ \frac12}_{loc}( \R_+, H^{-3}_{loc})$, 
we can use Corollary 7 of \cite{S} to get that $(A^\e)_{0 < \e <\e_0}$ is strongly compact in 
$L^2_{loc} (\R_+, L^2_{loc}) = L^2_{loc}(\R_+ \times \R)$. Since, by Lemma \ref{Bornes}, 
$A^\eps - u^\eps$ tends to zero strongly in $L^\infty(\R_+, L^2)$, we also get that 
$(u^\eps)_{0 < \e <\e_0}$ is strongly compact in $L^2_{loc}(\mathbb{R}_+, L^2_{loc})$.

\bigskip

Let now $A\in L^2_{loc}(\mathbb{R}_+, L^2_{loc})$ and $0 < \e_j \to 0$ as $ j\to +\ii$ such that
\begin{eqnarray}
& &\label{conv1}  A^{\e_j} \mbox{ converges to  }A \mbox{ strongly in }  L^2_{loc}(\R_+, L^2_{loc}), 
\mbox{ and weakly in } 
L^2_{loc}(\mathbb{R}_+, H^1_{loc});\\
& & \label{conv2} u^{\e_j} \mbox{ converges to  } A \mbox{ in }  L^2_{loc}(\R_+, L^2_{loc}).
\end{eqnarray}
Note that the weak convergence of $A^\eps$ just comes from the uniform $H^1$ bound which comes
 from Lemma \ref{Bornes}.

\bigskip

The next step in the proof is to obtain that  $A$ is a weak solution to the KdV equation.\\

For that purpose, let us write from \iref{PhAm} the equation satisfied by 
$A^{\e_j}+u^{\e_j}$ in the weak form:
\begin{eqnarray*}
& & \int_{\R_{+} \times \R} \big(
     A^{\e_j} + u^{\e_j} \big) \p_{t} \zeta \, dtdx + 
\int_{\R_{+} \times \mathbb{R}}\Big(
       (u^{\e_j})^2 + \frac{1}{\eps_j^4} g ( \eps_j^2 A^{\e_j}) 
- \frac{\eps_j^2 (\p_{x}A^{\eps_j} )^2}{4 c^2 (1 + \eps_j^2 A^{\eps_j})^2} \Big) \p_{x} \zeta\,dtdx \\
  & & + \int_{\R_{+} \times \R} \frac{ \p_{x}A^{\eps_j}}{4c^2 ( 1 + \eps_j^2 A^{\eps_j})} \p_{xx} \zeta\,dtdx
     +\int_{\R_{+}\times \R}\big(- u^{\eps_j} \p_{x} A^{\eps_j} \zeta  +
  A^{\eps_j} u^{\eps_j} \p_{x}\zeta\big)\,dtdx \\
 & &  = \int_{\mathbb{R}}\big(A_{0}^{\eps_j} + u_{0}^{\eps_j} \big) \zeta(0,x)\,dx
     \end{eqnarray*}
     for every $\zeta \in \mathcal{C}^\infty_{c}(\mathbb{R}\times \mathbb{R})$.
One can  pass to the limit easily  in most of the terms 
 by the strong convergence. Moreover, we can use that
$$ \int_{\R_{+}\times \R} u^{\eps_j} \partial_x A^{\eps_j} \zeta \to 
\int_{\R_{+}\times \R} A\, \p_{x} A \, \zeta $$
since $u^{\eps_j} \to A$ strongly and $\partial_{x} A^\eps \to \p_x A$ weakly. Since 
$ A^\eps$ is uniformly bounded in $L^\infty(\mathbb{R}_{+}, H^1)$, we have that
$$ \Big|\int_{\mathbb{R}_{+}\times \mathbb{R}} 
\frac{\eps_j^2 (\partial_{x}A^{\eps_j} )^2}{4 c^2 (1 + \eps_j^2 A^{\eps_j})^2} \p_{x} \zeta\,dtdx\Big| 
\leq K \eps_j^2 \rightarrow 0.$$
Moreover, since
$$ \int_{\R_{+} \times \R} \frac{ \p_{x} A^{\eps_j} } { 4c^2 ( 1 + \eps_j^2 A^{\eps_j})} \p_{xx} \zeta\,dtdx
 = \int_{\R_{+} \times \R} \frac{ \p_{x}A^{\e_j}}{4c^2 } \p_{xx} \zeta\,dtdx 
-\e_j^2 \int_{\R_{+} \times \R}\frac{ A^{\e_j} \p_{x}A^{\e_j} }{4c^2 (1 +\e_j^2 A^{\e_j})} \p_{xx} \zeta\,dtdx,$$
we get that the first term converges to
$$ \int_{\mathbb{R}_{+} \times \mathbb{R}} \frac{\p_{x} A}{4c^2 } \partial_{xx} \zeta\,dtdx$$
by weak convergence and that the second term converges to zero because of the uniform bounds. Therefore,
$$ \int_{\R_{+} \times \R}\frac{ \p_{x} A^{\eps_j}}{4c^2 ( 1 + \eps_j^2 A^{\eps_j})} \p_{xx} \zeta\,dtdx 
\ra \int_{\mathbb{R}_{+} \times \mathbb{R}} \frac{\p_{x} A}{ 4c^2 } \partial_{xx} \zeta\,dtdx. $$
Finally, we write
$$ \t{f}(r) = \big[ 1 + \frac{ 2}{c^2} f''(1) \big] r^2 + \BO(r^3) 
\quad \quad {\rm as} \quad r \to 0,$$
to infer
$$ \int_{\R_{+} \times \mathbb{R}} \frac{1}{\eps_j^4} \t{f}( \eps_j^2 A^{\e_j}) \p_{x} \zeta\,dtdx 
\to \big[ c^2 + 2 f''(1) \big] \int_{\R_{+} \times \mathbb{R}} A^2 \p_{x} \zeta\,dtdx $$

Consequently, we finally obtain that $A$ satisfies
$$  \int_{\R_{+}\times \R} \Big( 2 A \partial_{t}\zeta + \frac k2 A ^2\partial_{x}\zeta 
    +\frac{1}{4 c^2} { \partial_{x} A } \partial_{xx}\zeta\Big)\,dtdx
    = \int_{\mathbb{R}} 2A_{0}(x) \zeta(0,x)\, dx , $$
which is the weak form of the KdV equation.\\

\bigskip

 Next,  by passing to the limit in the bound 
of Lemma \ref{Bornes}, we get that $A \in L^\ii \big( \R_+, H^1 \big)$.
 Moreover, since it is a solution of the KdV equation, we deduce 
  that
  $$ \p_t A =  \frac{1}{8c^2} \p^3_x A - \frac{k}{2} A \p_x A \in L^\ii \big( \R_+, H^{-2} \big).$$
  Hence  $A \in$ Lip $(\R_+, H^{-2})$,  and  by interpolation in space, we get that 
$A \in \BC^0_b(\R_+,H^s)$ for any $0 \leq s < 1$.

We shall now prove that 
  $A=v$ the  unique solution 
 of the KdV equation given by  Theorem \ref{CauchyKdVH1}.
This fact  can be deduced from a general  uniqueness theorem for the KdV
 equation \cite{Zhou}. Nevertheless, here,   by using that
  the solution $v$ given by Theorem \cite{KPV}  verifies the additional property
   $\partial_{x}v \in L^4_{loc}(\mathbb{R}_{+}, L^\infty)$, one can get
    that $A=v$ by a very simple weak strong uniqueness argument.  Indeed, let us set
$\theta \equiv A-v$ and 
observe that $\theta \in L^\ii \big( \R_+, H^1 \big) \cap \BC^0_b(\R_+,H^s)$ for $0<s<1$ solves
$$ 2\p_t \theta - \frac{1}{4c^2} \p^3_x \theta = - k A \p_x \theta - k \theta \p_x v 
= - k \theta \p_x \theta -  k \theta \p_x v -k v \partial_{x}\theta , 
\quad \quad \quad \quad \theta_{|t=0}=0 .$$

Consequently, the standard $L^2$ energy estimate for this equation gives
\begin{align*}
\frac{d}{dt} \int_\R \theta^2 \ dx 
 \leq 2|k| |\!|\partial_{x} v |\!|_{L^\infty} \, |\!| \theta |\!|_{L^2}^2
\end{align*}
By the standard Gronwall inequality,  this yields immediately that $\theta =0$, since $\theta _{|t=0}=0$,  
and $\partial_{x}v \in L^4_{loc}( L^\infty)\subset L^1_{loc}(L^\infty).$

\bigskip

As a consequence of the uniqueness of the limit, the full sequence $A^\e $ converges to $v$ 
as $\e \to 0$ strongly in $L^2_{loc}(\R_+, L^2_{loc})$ and weakly in $L^2_{loc}(\R_+, H^1_{loc})$, 
where $v$ is {\it the} $H^1$-solution of the KdV equation of Theorem \ref{CauchyKdVH1}.\\

It remains to improve the convergence of $A^\eps$ i.e. to prove that 
 we actually have   the local in time global in space strong convergence, as $\e \to 0$,
$$ A^\e \to v \quad \quad {\rm in} \quad \BC \big( [0, T], L^2 \big)$$
for every $T>0$.

\bigskip

From Lemma \ref{Bornes} and the proof of Lemmas \ref{init} and \ref{source}, we infer that
$$ A^\e + u^\e \quad \mbox{is uniformly bounded in} \ 
L^\ii \big( \R_+, L^2 \big) \cap {\rm Lip} \big( \R_+ , H^{-2} \big). $$
In particular,
$$ A^\e + u^\e \quad \mbox{is uniformly bounded in} \ 
\BC^{0,1/2} \big( \R_+, H^{-1} \big) \cap L^\ii \big( \R_+, L^2 \big). $$
Since we already have that
$$
A^\e + u^\e \to 2A=2v \quad \quad {\rm in} \quad L^2_{loc}\big(\R_+, L^2_{loc}\big) ,
$$
it follows by a new use of the Aubin Lions lemma  that
\be
\label{nv}
 A^\e + u^\e \to 2v \quad \quad {\rm in} \quad \BC^0_{loc} \big( \R_+, H^{-1}_{loc} \big).
 \ee
Consequently,  we can write for  every $T>0$, $R>0$, 
$$ \sup_{[0, T]} |\!|A^\eps - v |\!|_{H^{-1}(-R,R)}
 \leq \frac{ 1}{ 2} \sup_{[0, T]}\Big(  |\!| A^\eps+ u^\eps - 2 v |\!|_{H^{-1}(-R, R) }
  +   |\!| A^\eps -  u^\eps |\!|_{H^{-1}(-R, R) }\Big)$$
  and since by Lemma \ref{Bornes}, we have that  $A^\e - u^\e \to 0$ in $L^\ii(\R_+,L^2)$,
  this yields  thanks to \eqref{nv} that 
$$ A^\e \to v \quad \quad {\rm in} \quad \BC^0_{loc} \big( \R_+, H^{-1}_{loc} \big). $$
Let us now  fix $T> 0$. We then prove that, as $\e \to 0$,
$$ \sup_{[0,T]} \big| \langle A^\e - v, v \rangle_{L^2} \big| \rightarrow 0.$$
Indeed, let $\eta >0$ be given. Since $v \in \BC^0_b \big( \R_+, L^2 \big)$, there 
exists $R > 0$ such that
$$ \sup_{[0,T]} \int_{|x| \geq R} v^2\ dx \leq \eta^2.$$
Next,   with  $\zeta \in \mathcal{C}^\infty_{c}( -2R, 2R)$ such that $\zeta=1
$ on $[-R, R]$, we split
$$ \sup_{[0,T]} \big| \langle A^\e - v,  v \rangle_{L^2} \big| \leq 
\sup_{[0,T]} \big| \langle A^\e - v, \zeta v \rangle_{L^2} \big| 
+ \sup_{[0,T]} \big| \langle A^\e - v, (1- \zeta) v \rangle_{L^2} \big|.$$
The first term tends to $0$ as $\e \to 0$ since $\zeta v\in \BC_b \big( \R_+, H^1 \big)$
 is compactly supported  and 
$A^\e \to v $ in $\BC^0_{loc} \big( \R_+, H^{-1}_{loc} \big)$. The second term is 
$\leq \eta \sup_{[0,T]} \big|\!\big| A^\e - v \big|\!\big|_{L^2} \leq K \eta $, and 
the limit follows.

Therefore,
\be
\label{no1} \sup_{[0,T]} \big|\!\big| A^\e - v \big|\!\big|_{L^2}^2 = 
\sup_{[0,T]} \Big\{ 
\big|\!\big| A^\e \big|\!\big|_{L^2}^2 - \big|\!\big| v \big|\!\big|_{L^2}^2 
- 2 \langle A^\e - v, v \rangle_{L^2} \Big\} = \sup_{[0,T]} \Big\{ 
\big|\!\big| A^\e \big|\!\big|_{L^2}^2 - \big|\!\big| v \big|\!\big|_{L^2}^2 \Big\}  + o(1).\ee
We now use that $E^\e(\psi^\e)$ and $\BI_0(A) = \int_\R A^2 \ dx$ are independent of $t$, thus
\be
\label{no2}
 \big|\!\big| A(t) \big|\!\big|_{L^2}^2 = \big|\!\big| A_0 \big|\!\big|_{L^2}^2 \ee
and,  using Lemma \ref{Bornes}  and the same expansion  as in Step 1 of  the proof of Lemma \ref{Bornes}, we infer
$$ E^\e\big( \psi^\e(t) \big) = 
\frac{\e^2}{2} \int_\R 4 c^2 \big( A^\e(t) \big)^2 + \big( \p_x \vp^\e(t) \big)^2 \ dx + \BO( \e^4) 
= 4 c^2 \e^2 \int_\R \big( A^\e(t) \big)^2 \ dx + \BO(\e^3).$$
 Note that  the $\BO(\e^3)$ is uniform with 
respect to $t\in \R_+$. Since $ E^\e\big( \psi^\e(t) \big) = E^\e\big( \psi^\e_0\big) $ and 
the same expansion holds at $t=0$ (this is Step 1 in the proof of Lemma \ref{Bornes}), we deduce
\be
\label{no3}
 \int_\R \big( A^\e(t) \big)^2 \ dx = \int_\R \big( A^\e_0 \big)^2 \ dx + \BO(\e),
 \ee
where $\BO(\e)$ is uniform with respect to $t\in \R_+$. Consequently, thanks to \eqref{no1}, \eqref{no2},
 \eqref{no3}, we obtain that 
$$ \sup_{[0,T]} \big|\!\big| A^\e - v \big|\!\big|_{L^2}^2 = 
\big|\!\big| A^\e_0 \big|\!\big|_{L^2}^2 - \big|\!\big| A_0 \big|\!\big|_{L^2}^2 + o(1),$$
and since $A^\e_0 \to A_0$ in $L^2$ by assumption, the result  in $L^2$ follows.\\

The proof of Theorem \ref{KdVH1} is now complete, since the convergence of $A^\eps$  
in $L^\ii_{loc}(\R_+,H^s)$, $0 < s < 1$ follows by interpolation in space using the convergence in 
$L^\ii_{loc}(\R_+,L^2)$ and the uniform  bounds in $L^\ii(\R_+,H^1)$. 

\bigskip

\subsection{Convergence in $\bs{H^1}$}

\label{H1+}

\ \indent In this subsection, we  shall  put a more restrictive assumption on the initial
 data, namely 
$$\big|\!\big| \partial_{x} \varphi_{0}^\eps - 2cA_{0}^\eps \big|\!\big|_{L^2}= o(\eps)$$
instead of $\BO(\e)$ in order to get the strong convergence in $H^1$ of the amplitude
 $A^\eps $. 

\begin{theo}
\label{conver}
Under the assumptions of Theorem \ref{KdVH1}, if,   at the initial time, we have the 
additional assumptions
$$ A_{0}^\eps \to A_0 \quad \quad \quad {\it in} \ \ H^1$$
and
 \be
 \label{boncondi} 
\big|\!\big| \partial_{x} \varphi_{0}^\eps - 2cA_{0}^\eps \big|\!\big|_{L^2}= o(\eps), \ee
then
$$ A^\e \to A \quad \quad {\rm in} \quad \BC^0_{loc} \big( \R_+,H^1(\R) \big).$$
\end{theo}
\bigskip

\subsection*{Proof.} 

The idea follows the one in the end of the proof of Theorem \ref{KdVH1}, but relies on the 
conservation of 
$$\BI_1\big( A(t) \big) \equiv \int_\R \frac{1}{ 4 c^2} (\partial_{x} A)^2 + \frac{k}{3 } A^3 \ dx$$
 for KdV and $ E_\e\big( \psi^\e(t) \big) - 2 c P^\e \big( \psi^\e(t) \big)$
 for \eqref{NLSd}. First, we expand to third order
$$ F(R) = c^2 (R-1)^2 + \frac{1}{3} f''(1)(R-1)^3 + F_4(R), \quad {\rm with} \quad F_4(1+r)=\BO(r^4), 
\quad r \to 0, $$
so that \iref{develoE-cP} becomes now
$$ E^\e(\psi) - 2c P^\e(\psi) = \frac{1}{2} \int_\R 
\big( \rho^2 -1 \big) (\p_x \phi)^2 + (\p_x \rho)^2 + 
\Big( \p_x \phi - \frac{c}{\e} (\rho^2-1 ) \Big)^2 + \frac{f''(1)}{3\e^2}\, \big( \rho^2-1 \big)^3 
+ \frac{1}{\e^2}\, F_4 \big( \rho^2-1 \big) \ dx. $$
Since $\p_x \vp^\e_0 - 2 c A^\e_0 = \BO(\e)$ in $L^2$  by assumption, we infer
$$ \int_\R \Big( \p_x \vp^\e_0 - 2 c A^\e_0 - c \e^2 ( A^\e_0 )^2 \Big)^2 \ dx = 
\int_\R \Big( \p_x \vp^\e_0 - 2 c A^\e_0 \Big)^2 \ dx + \BO(\e^3) = o(\e^2),$$
using the hypothesis \iref{boncondi}. Therefore, at time $t=0$, we infer, as in Step 1 of the 
proof of Lemma \ref{Bornes}, that
\begin{align*}
E^\e ( \psi^\e_0 ) - 2c P^\e(\psi^\e_0) & \ = \frac{\e^4}{2} \int_\R 
2 A^\e_0 \big( \p_x \vp^\e_0 \big)^2 + ( \p_x A^\e_0 )^2 + \frac{8 f''(1)}{3} ( A^\e_0)^3 \ dx \\ 
& \quad \quad 
+ \frac{\e^2}{2} \int_\R \Big( \p_x \vp^\e_{0} - 2 c A^\e_{0} - c \e^2 ( A^\e_{0} )^2 \Big)^2 \ dx 
+ \BO(\e^6) \\ & \ = \frac{\e^4}{2} \int_\R ( \p_x A^\e_0 )^2 + 
8 \Big[ c^2 + \frac{f''(1)}{3} \Big] ( A^\e_0)^3 \ dx + o(\e^4)
\\ & \ = 2c^2 \e^4 \BI_1 \big( A^\e_0 \big) + o(\e^4) = 2c^2 \e^4 \BI_1 \big( A_0 \big) + o(\e^4),
\end{align*} 
since $A^\e_0 \to A_0$ in $H^1(\R) \subset L^3(\R)$. Similarly, 
given $t\in \R_+$ and using Lemma \ref{Bornes}, we have
$$ E^\e\big( \psi^\e(t) \big) - 2c P^\e \big( \psi^\e(t) \big) = 
2c^2 \e^4 \BI_1 \big( A^\e(t) \big) + 
\frac{\e^2}{2} \int_\R \Big( \p_x \vp^\e(t) - 2 c A^\e(t) \Big)^2 \ dx + \BO(\e^5),$$
where $\BO(\e^5)$ is uniform with respect to time. Since $\BI_1 \big( A(t) \big)$ and 
$E^\e(\psi^\e) - 2c P^\e(\psi^\e)$ are independent of time, this implies,
\be
\label{I1conserve}
\BI_1 \big( A(t) \big) = \BI_1 \big( A^\e(t) \big) + 
\frac{1}{4c^2\e^2} \int_\R \Big( \p_x \vp^\e(t) - 2 c A^\e(t) \Big)^2 \ dx + o(1)
\ee
uniformly in time. 

Now, let us study  the term involving the 
 $L^3$-norm in $\BI_1$. Let $T>0$ be fixed. From Lemma \ref{Bornes}, 
$A^\e$ is uniformly bounded in $L^\ii(\mathbb{R}_{+} \times \mathbb{R})$. Moreover, we have 
proved in Step 4 that $A^\e \to A$ in $\BC\big( [0,T], L^2 \big)$. As a consequence, 
$A^\e \to A$ in $\BC \big( [0,T],L^3 \big)$. 
Inserting this in \iref{I1conserve} yields, uniformly for $t\in [0,T]$,
\be
\label{I1conserve2}
\int_\R \big( \p_x A(t) \big)^2 \ dx = \int_\R \big( \p_x A^\e(t) \big)^2 \ dx + 
\frac{1}{\e^2} \int_\R \Big( \p_x \vp^\e(t) - 2 c A^\e(t) \Big)^2 \ dx + o(1).
\ee
We now consider
$$ \nu^\e (T) \equiv \sup_{[0,T]} \Big\{ 
\big|\!\big| \p_x A^\e - \p_x A \big|\!\big|_{L^2}^2 + 
\frac{1}{\e^2} \big|\!\big|\p_x \vp^\e - 2 c A^\e \big|\!\big|_{L^2}^2 \Big\}.$$
Since $A \in \BC([0,T],H^1)$, arguing as in the end of the proof of Theorem \ref{KdVH1}, 
we infer
$$ \nu^\e(T) = \sup_{[0,T]} \Big\{ 
\big|\!\big| \p_x A^\e \big|\!\big|_{L^2}^2 - \big|\!\big| \p_x A \big|\!\big|_{L^2}^2 + 
\frac{1}{\e^2} \big|\!\big|\p_x \vp^\e - 2 c A^\e \big|\!\big|_{L^2}^2 \Big\} + o(1).$$
Combining this with \iref{I1conserve2} gives $\nu^\e(T) = o(1)$ as desired. 
This ends the proof of Theorem \ref{conver}.

\section{The general $n$ dimensional case}

\subsection{Proof of Theorem \ref{Smoothy}}

\label{section2}

\ \indent It is more convenient  to use a different hydrodynamic form of (NLS).
 As in \cite{Grenier}, we shall seek for a solution of \eqref{NLSd} under
  the form
  \be
  \label{grenier}
\psi^\eps = \big( 1+ \eps^2 a^\e (t,X) \big) e^{i \eps \theta^\e(t,X)}, 
\quad a^\e \in \mathbb{C}, \quad \theta^\e \in \R, \quad \eps^2 |a^\e| \leq \frac{ 1}{2}\ee
  that is to say that we allow the amplitude to be complex at positive times.
    The reason for this choice is that we can obtain an hydrodynamic equation 
     for $(a^\eps, \theta^\eps)$ which is much simpler.
   We shall prove that $a^\eps$ and $\theta^\eps$ are well defined on $[0,T]$ for some 
$T>0$ independent
    of $\eps$ and satisfy for $s>1+n/2$  the uniform estimate
  \be
  \label{unif2}
 \big|\!\big| a^\eps(t) \big|\!\big|_{H^{s+1}} + \big|\!\big| \p_x \theta^\eps(t) \big|\!\big|_{H^{s+1}} 
+ \e \big|\!\big|\nabla_\perp \theta^\eps(t) \big|\!\big|_{H^{s+1}} \leq C, 
\quad \forall t \in[0, T], \, \forall \eps \in (0,    \eps_{0}]
   \ee
   for some $C>0$ independent of $\eps$.

   Note that once this estimate is proven, the  representation \eqref{polaire} and the estimate
    \eqref{unif} immediately follow.
    Indeed, for $\eps$ sufficently small, we get that $|\psi^\eps|$ remains far from zero
     on $[0, T]$ and  we have the relations
   \be
   \label{relat}
 A^\eps = \frac{ |1+ \eps^2 a^\eps| - 1 }{ \eps^2 }, \quad \quad \p_j \varphi^\eps= 
    \p_j \theta^\eps + \frac{\eps}{i} \Big( \frac{\p_j a^\eps}{ 1+ \eps^2 a^\eps}
     - \frac{ \p_j A^\eps}{1+ \eps^2 A^\eps} \Big) \quad 1\leq j\leq n \ee
     from which we deduce by standard  Sobolev-Gagliardo-Nirenberg-Moser estimates that
$$ \big|\!\big|A^\eps(t)\big|\!\big|_{H^{s+1}} + 
\big|\!\big| \p_x \varphi^\eps(t) \big|\!\big|_{H^{s}} + 
\e \big|\!\big| \nabla_\perp \varphi^\eps(t) \big|\!\big|_{H^{s}} \leq C\, \quad \forall t \in [0,T], 
     \quad  \forall \eps \in  (0, \eps_{0}]$$
      for some $C$ independent of $\eps$ since $s>1+n/2$.
      
  Let us now write down  the equation for $(a^\eps, \theta^\eps)$.
   By plugging the anzatz \eqref{grenier} in \eqref{NLSd}, we get
   \begin{eqnarray*}
   & &  i c \eps^3 \Big( \eps^2 \partial_{t} a^\eps  + i \eps( 1 + \eps^2 a^\eps ) \partial_{t} \theta^\eps
    \Big) - i c \eps \Big(  \eps^2 \partial_{x} a^\eps +i  \eps ( 1 + \eps^2 a^\eps  ) \partial_{x} \theta^\eps \Big) \\
    & &  + \frac{ \eps^2}{ 2 } \Big( \eps^2 \Delta^\eps a^\eps + 2 i \eps^3 \nabla^\eps \theta^\eps \cdot \nabla^\eps a^\eps
     + i \eps ( 1 + \eps^2 a^\eps ) \Delta^\eps \theta^\eps - \eps^2 ( 1+ \eps^2 a^\eps )|\nabla^\eps \theta^\eps |^2
     \Big) \\
      & &  - ( 1 + \eps^2 a^\eps ) f\big( |1+ \eps^2 a^\eps |^2 \big) = 0
       \end{eqnarray*}
       where we use the notation
       $$ \nabla^\eps = (\partial_{x}, \eps \partial_{\perp})^t, \quad \Delta^\eps=  
\nabla^\eps \cdot \nabla^\eps  = \partial_{x}^2 + \eps^2 \Delta_{\perp}.$$
      Since  we allow the amplitude $a^\eps$ to be complex, we have some freedom
       to write down  hydrodynamic equations. As noticed  in \cite{Grenier}, it is convenient to split the above
        equation into         the system
     \begin{eqnarray*}
    \left\{ \begin{array}{ll}
     & \displaystyle{  \partial_{t} a^\eps - \frac{ 1}{\eps^2} \partial_{x}a^\eps  + \frac{ 1}{c} \nabla^\eps \theta^\eps
      \cdot \nabla^\eps a^\eps + \frac{ 1}{ 2 c \eps^2} ( 1 + \eps^2 a^\eps ) \Delta^\eps \theta^\eps
       = \frac{i}{ 2 \eps c} \Delta^\eps a^\eps} \\       \\
       & \displaystyle{ \partial_{t} \theta^\eps - \frac{ 1}{ \eps^2 } \partial_{x} \theta^\eps  + 
\frac{ 1 }{ 2 c} |\nabla^\eps  \theta^\eps |^2 + \frac{ 1}{ c\eps^4} f\big( |1+ \eps^2 a^\eps |^2 \big)}=0.
      \end{array}\right.
       \end{eqnarray*}
Consequently, by using the new unknown $v^\eps \equiv \ds{\frac{ 1}{2 c} \nabla^\eps \theta^\eps}$, 
we get
    \be
     \label{sgrenier}
    \left\{ \begin{array}{ll} 
    &  \displaystyle{  \partial_{t} a^\eps - \frac{ 1}{ \eps^2} \partial_{x}a^\eps  +  2 v^\eps
      \cdot \nabla^\eps a^\eps + \frac{1}{ \eps^2} ( 1 + \eps^2 a^\eps ) \nabla^\eps\cdot v^\eps  
       = \frac{i}{ 2 \eps c} \Delta^\eps a^\eps} \\ \\
     & \displaystyle{ \partial_{t} v^\eps  - \frac{ 1}{\eps^2 } \partial_{x} v^\eps  + 
2 v^\eps \cdot \nabla^\e v^\eps  + \frac{ 1}{ 2 c^2 \eps^2} f'\big( |1+ \eps^2 a^\eps |^2 \big) \Big( 2 \nabla^\eps \mbox{Re }a^\eps+ \eps^2 \nabla^\eps |a^\eps |^2}\Big)=0.
      \end{array}\right. 
      \ee
      We add to this system the initial condition
      \be
      \label{initgrenier}
      a^\eps(0, X)= A_{0}^\eps(X), \quad v^\eps(0, X)= \frac{ 1}{ 2c }\nabla^\e \vp_{0}^\e (X).
      \ee
    Consequently, 
     we can set
      $U^\eps\equiv (\mbox{Re }a^\eps,\, \mbox{Im }a^\eps,  v^\eps)^t\in \mathbb{R}^{2+n}$, $\partial^\eps=
      (\partial_{x}, \eps \partial_{\perp})$
       and write the above system under the abstract form:
    \be
   \label{abstr}
   \partial_{t}U^\eps + \frac{ 1}{\eps^2} H( \eps^2 U^\eps, \partial^\eps) U^\eps =  \frac{ 1}{ \eps} L(\partial^\eps) U^\eps
   \ee
     where $L(\partial^\eps)$ is a constant coefficients  second  order differential operator
    $$ L(\partial^\eps) \equiv \frac{ 1}{ 2 c } \left( \begin{array}{ccc}  J \Delta^\eps  &  0  \\
     0     &   0 
     \end{array}\right), \quad J = \left(  \begin{array}{cc} 0 & -1 \\ 1 &  0 \end{array}\right)
   $$     
     and $H(\eps^2 U^\eps, \partial^\eps)$ is a first order hyperbolic operator
   $$  H(\eps^2 U^\eps, \partial^\eps) = \sum_{k=1}^n H^k(\eps^2 U^\eps) \partial_{k}^\eps, \quad $$
    with symbol
    $$ H(\eps^2 U^\eps, \xi)= \sum_{k=1}^n H^k(\eps^2 U^\eps) \xi_{k} = 
    \left( \begin{array}{cc}
     ( -\xi_{1}+ 2  \eps^2 v^\eps \cdot \xi)I_{2} &  ( e + \eps^2 a^\eps)  \xi^t \\
             \big (1 + g( \eps^2  a^\eps) \big)\xi\big( e + \eps^2 a^\eps \big)^t & 
\big( -\xi_{1}+ 2  \eps^2 v^\eps \cdot \xi \big) I_{n} \end{array}\right)$$
        where 
        $$ e \equiv \left(  \begin{array}{ll} 1 \\ 0 \end{array} \right)$$
        and $g$ is defined by the expansion:
\be
\label{expd}
 \frac{1}{ c^2}f'(1+  2(a, 1)  +  |a|^2)= 1 + g(a), \quad g(a)= \mathcal{O}(|a|), \, |a|\leq 1
 \ee 
 since $f'(1) = c^2$.
 
 Note that the structure of \eqref{abstr} is much simpler than the one of the  standard
  hydrodynamic system for $( A^\e, \nabla^\eps \varphi^\eps)^t$ that is obtained 
from \eqref{PhAmd} by the standard Madelung transform. Indeed, \eqref{abstr} is a simple skew-symmetric  constant coefficient perturbation of an hyperbolic system.

  Note that the difficulties du to the presence of vacuum
     which arise in the study of NLS with solutions which tends to zero
      at infinity (\cite{Alazard-Carles}, \cite{Chiron-Rousset} are not present here. The above system can be  easily symmetrized by using
   $$ S(\eps^2 U^\eps) = \left( \begin{array}{ccc}  I_{2 }& 0 
    \\ 0  &  \frac{ 1 }{ 1 + g(\eps^2 a^\eps) } \end{array} \right) $$
    which is positive. 
   Indeed, we have
$$ S(\eps^2 U^\eps) L(\partial^\eps) = \frac{ 1}{ 2 c }\left( \begin{array}{ccc}J \Delta^\eps & 0 \\ 
0 & 0 \end{array} \right) $$ 
    which is a skew symmetric operator:
    \be
    \label{sym1}
    \Big( S^\eps(\eps^2 U^\eps) L(\partial^\eps)  V, V \Big)=0, \quad \forall V \in H^2(\mathbb{R}^n)
    \ee
    where we use the notation $(\cdot, \cdot)$ for the $L^2(\mathbb{R}^n)$ scalar product.
    Moreover, we also have that
$$ S(\eps^2 U^\eps) H(\eps^2U^\eps, \xi) = 
    \left( \begin{array}{cc}
     ( -\xi_{1}+ 2  \eps^2 v^\eps \cdot \xi)I_{2} &  ( e + \eps^2 a^\eps)  \xi^t \\
             \xi \big( e + \eps^2 a^\eps \big)^t & 
\ds{\frac{ 1}{1 + g(\eps^2 a^\eps ) }}\big( -\xi_{1}+ 2 \eps^2 v^\eps \cdot \xi
        \big) I_{n} \end{array}\right) $$
   is symmetric for every $\xi \in \mathbb{R}.$

The local existence and uniqueness  of a smooth solution $U^\eps \in \mathcal{C}([0, T^\eps),  H^{s+1})$
for this system  is classical.
 Moreover,  let us define
 $$ T^\eps_{*}= \sup \Big\{ T \in [0, T^\eps), \, \, \forall t \in [0, T], \quad 
  |\eps^2 a^\eps  |_{L^\infty} \leq \frac{1}{ 2}, \quad  |\!|U^\eps |\!|_{ H^{s+1}}<+\infty \Big\}.$$ 
 We shall   prove that  $T^\eps_{*}$ is bounded from below by a positive number  when $\eps$ tends to 
 zero. This will be achieved by proving  $H^{s+1}$ estimates uniform  in $\eps$.
 
 Note that for $t \leq T^\eps_{*}$, the symmetrizer $S(\eps^2 U^\eps)$ is well defined
  and verifies 
\be
   \label{sym3}
   \big( S( \eps^2 U^\eps) V, V \big) \geq c_{0} |\!|V|\!|_{L^2}^2 , \quad \forall t \in [0, T^\eps_{*}],
    \quad \forall V\in L^2(\mathbb{R}^n)
   \ee    
 for some $c_{0}>0$ independent of $\eps$.  
       Moreover, thanks to an integration by parts, we also   have for some $C>0$ independent of $\eps$
     that 
    \be
    \label{sym2}
   \big|  \big(S(\eps^2 U^\eps) H(\eps^2 U^\eps, \partial^\eps)V, V \big) \big|
    \leq C \eps^2 \big|\!\big| \nabla U^\eps \big|\!\big|_{L^\infty} \big|\!\big| V\big|\!\big|_{L^2}^2, \quad
     \forall t \in [0, T^\eps_{*}]
   \ee
   for every $V\in H^1(\mathbb{R}^n)$.   
   
    We can  now easily perform  for $s>1+ n/2$ an $H^{s+1}$  estimate  for \eqref{abstr}.
  Indeed, for every $\alpha \in \mathbb{N}^n$, $|\alpha | \leq s+1$, we have
  \be
  \label{abstr2}
 \partial_{t}\partial^\alpha U^\eps + \frac{ 1}{ \eps^2} H( \eps^2 U^\eps, \partial^\eps) \partial^\alpha U^\eps 
- \frac{ 1}{ \eps} L(\partial^\eps) \partial^\alpha U^\eps 
+ \frac{ 1}{\eps^2} \big[ \partial^\alpha, H(\eps^2 U^\eps,\partial^\eps ) \big] U^\eps =0
  \ee
  By the standard tame Gagliardo-Nirenberg-Moser estimate, we get that
  \be
  \label{sym4}
  \Big|\!\Big| \frac{1}{ \eps^2} \big[ \partial^\alpha, H(\eps^2 U^\eps,
  \partial^\eps ) \big] U^\eps \Big|\!\Big|_{L^2} 
\leq C \big|\!\big| U^\eps \big|\!\big|_{W^{ 1, \infty}} \big|\!\big|U^\eps \big|\!\big|_{H^{s+1}},
   \quad \forall t \in [0, T^\eps_{*}].
  \ee
  From now on $C$ is a number independent of $\eps$ which may change from line to line.
   
  By using \eqref{sym1}, \eqref{sym2}, \eqref{sym4}, we get the energy estimate: 
  $$ \frac{ d}{ dt } \Big( \frac{ 1}{2 } \big( S(\eps^2 U^\eps ) \partial^\alpha U^\eps, 
   \partial^\alpha U^\eps)  \Big)  \leq C \Big( \eps^2  |\!| \partial_{t} U^\eps |\!|_{L^\infty}
    + \big|\!\big| U^\eps \big|\!\big|_{W^{1, \infty}} \Big)\big |\!\big| U^\eps \big|\!\big|_{H^{s+1}}^2, 
\quad \forall t \in [0, T^\eps_{*}] .$$
   By  using  \eqref{abstr}, we get that
   $$ \big|\!\big|\partial_{t} U^\eps \big|\!\big|_{L^\infty} \leq 
C \Big( \frac{1}{\eps^2} \big|\!\big|U^\eps \big|\!\big|_{W^{1, \infty}}
    + \frac{ 1}{ \eps } \big|\!\big|U^\eps \big|\!\big|_{W^{2, \infty}}\Big).$$ 
    Consequently, we can   integrate in time and use \eqref{sym3} to get
 \be
 \label{abstr3}
 \big|\!\big| U^\eps(t) \big|\!\big|_{H^{s+1}}^2 \leq C \Big( 
\big|\!\big| U^\eps_{0} \big|\!\big|_{H^{s+ 1 }}^2 + \int_{0}^t 
  \big|\!\big| U^\eps \big|\!\big|_{W^{ 2, \infty}} 
\big|\!\big|U^\eps(\tau) \big|\!\big|_{H^{s+ 1 }}^2\, d\tau\Big).
  \ee
  Finally, by using the  Sobolev embedding $ H^{s+1} \subset  W^{2, \infty}$  for
   $s>1+ n/2$, we find in a classical way from \eqref{abstr3} that $T^\eps_{*}>T>0$ 
for every $\eps \in  (0, \eps_{0})$ for some $\eps_{0}$ 
sufficiently small. We refer for example to \cite{Klainerman-Majda}, \cite{GrenierS}, \cite{Schochet}
 for more details. This ends 
 the proof of Theorem \ref{Smoothy}.
    
  \bigskip 
 
 \subsection{Proof of Theorem \ref{asympt1}}  
 \label{sectionas}
    We shall now study  the convergence towards the KP-I  equation. We could pass 
to the limit directly from \eqref{sgrenier}. Nevertheless, to make a link more clear 
with the first part of the paper and the formal derivation, we shall  pass to the limit 
directly from the standard hydrodynamic equation \eqref{PhAmd}.
      As already explained in  the beginning of the proof, we can deduce
       from the representation \eqref{grenier} and the bounds \eqref{unif2}
        that the smooth representation \eqref{polaire} with the uniform bounds \eqref{unif}
         hold on $[0, T]$.
         Consequently, we already have 
    \be
    \label{Au}
    \big|\!\big|  A^\eps(t) \big|\!\big|_{H^{s+1 }} + \big|\!\big|u^\eps(t) \big|\!\big|_{H^{s}} 
\leq C, \quad \forall t \in [0, T], \, \forall \eps \in (0,\eps_{0} )
    \ee
    for $s>1+n/2$, where $(A^\eps, u^\eps = \frac{ 1 }{2 c} \nabla^\eps \varphi^\eps)$ 
solves the system
    \be
\label{euler1}
\left\{
\begin{array}{ll}
\displaystyle{ \p_{t }A^\eps - \frac{ 1}{ \eps^2 }\partial_{x}A^\eps  
+ \frac{ 1}{ \eps^2} \nabla^\eps \cdot u^\eps
 +  2 u^\eps \cdot \nabla^\eps  A^\eps + A^\eps \nabla^\eps  \cdot u^\eps = 0 } \\
 \\
  \displaystyle{ \partial_{t}u^\eps - \frac{ 1}{ \eps^2} \partial_{x}u^\eps +
   \frac{ 1}{ \eps^2}  \nabla^\eps A^\eps + 2  u^\eps \cdot \nabla^\eps  u^\eps + 
\frac{ 1}{\eps^2}g(\eps^2 A^\eps) \nabla^\eps A^\eps  = \frac{ 1}{ 4 c^2}\nabla^\eps 
\Big( \frac{\Delta^\eps A^\eps }{  1 + \eps^2 A^\eps } \Big)}.
   \end{array}\right.
   \ee
    Note that 
     $ \nabla^\eps \times u^\eps=0$,  hence, we obtain in particular that
     \be
     \label{cross}
     \partial_{x} u_{\perp}^\eps  = \eps\nabla_\perp u_{1}^\eps.
     \ee

   We can  apply  $\partial_{x}$ to the first equation and the first line of the second
      equation in \eqref{euler1} to get the system:
    \be
    \label{transsmooth}
 \left\{  \begin{array}{ll} \displaystyle{ \partial_{t} \partial_{x} A^\eps 
+ \frac{ 1}{\eps^2} \partial_{x}\big( \partial_{x} u^\eps_{1} - \p_{x} A^\eps \big) = {\rm S}_{A}^\eps} \\  \\
\displaystyle{ \partial_{t} \partial_{x} u_{1}^\eps 
+ \frac{ 1}{\eps^2} \partial_{x} \big( \partial_{x} A^\eps - \partial_{x} u^\eps_{1} \big) = {\rm S}_{u}^\eps,} \end{array}\right.
 \ee
    where 
    \begin{eqnarray*}
  & & {\rm S}_{A}^\eps \equiv - \partial_{x}\big( 2 u^\eps \cdot \nabla^\eps A^\eps 
+  A^\eps \nabla^\eps \cdot u^\eps \big) - \frac{ 1}{\eps} \partial_{x} \nabla_\perp \cdot u^\eps_{\perp} \\
& & {\rm S}_{u}^\eps \equiv - \partial_{x}\Big( 2 u^\eps \cdot \nabla^\eps u^\eps_{1} + \frac{ 1}{\eps^2} g(\eps^2 A^\eps) \partial_{1}
  A^\eps\Big)  +\frac{ 1}{4c^2} \partial_{xx} \Big( \frac{ \Delta^\eps A^\eps}{ 1 + \eps^2 A^\eps} \Big).
  \end{eqnarray*}
  By using \eqref{cross} and the $H^{s+1}$ bound \eqref{Au} which holds for $s>1+n/2 \geq 3/2$,
      we get the uniform estimate
$$ \big|\!\big|( {\rm S}_{A}^\eps ,{\rm S}_{u}^\eps) \big|\!\big|_{H^{-2}} \leq C, \quad \forall 
t \in[0, T], \quad \forall \eps \in (0, \eps_{0}] $$
     for some $C>0$.
     
   Consequently, from the proof of Lemma \ref{source} (it suffices to integrate also  with respect
    to the transverse variable), we get that:
     $ \partial_{x} A^\eps$ and $\partial_{x} u^\eps_{1}$ are uniformly bounded
      in $H^{\frac{1}{2}}( 0, T, H^{-3}_{loc})$ and also (see  \eqref{L2}) that
      \be
      \label{illp}
\partial_{x}A^\eps -\partial_{x}u^\eps_{1}= \mathcal{O}(\eps) \quad \mbox{  in } L^2(0, T, H^{-2}_{loc}).
      \ee
       Consequently, we can use again the relative compactness criterion of  \cite{S} and \eqref{Au}
        to get that $\partial_{x}A^\eps$ is strongly compact in 
         $L^2(0, T, H^m_{loc})$ and $\partial_{x}u^\eps_{1}$ in $L^2(0, T, H^{m-1}_{loc})$ for every $m<s$.
          Note that since $s>1$, one can  choose
          $ m>1$.  Consequently, 
          the way to recover the weak form of the KP-I  or KdV equation will be very close
           to what was done in the proof of Theorem \ref{KdVH1}. 
          We can take a subsequence $\eps_{j} \to 0$ such that 
\begin{eqnarray*}
& & \partial_{x}A^{\eps_j} \rightarrow \partial_{x }A \quad \mbox{ strongly in } L^2\big(0, T, H^m_{loc}\big), 
\quad \quad  \partial_{x}u^{\eps_j}_{1} \rightarrow \partial_{x }u_{1}\quad  \mbox{strongly  in } L^2\big(0, T, H^{m-1}_{loc}\big), \\
& & \ \ \ A^{\eps_{j}}  \rightarrow A \quad \ \ \, \mbox{ weakly in }  L^2\big(0, T, H^{s+1}\big), 
\hspace{1.25cm} u^{\eps_{j}}  \rightarrow  u \quad \quad \mbox{ weakly in }  L^2\big(0, T, H^{s}\big)
          \end{eqnarray*}  
        and moreover,  from \eqref{illp}, we also have
       \beq
       \label{limprop}
        A=  u_{1} \quad \quad  \mbox{ for almost every } t \in [0, T], \,X \in \mathbb{R}^n.
        \eeq
      As in the proof of  Theorem \ref{KdVH1}, the above properties are sufficient
       to pass to the limit in the weak form of the equation satisfied
        by $\partial_{x}A^\eps + \partial_{x} u^\eps_{1}$. Indeed,  by using \eqref{cross}, 
         we get from \eqref{euler1} that 
        \begin{eqnarray*}
    & & \int_{[0, T] \times \mathbb{R}^n}
     \partial_{x}\big(A^{\eps_{j}} + u_{1}^{\eps_{j}} \big)
      \partial_{t} \zeta + \Big(   2 u_{1}^{\eps_{j}}\big(  \partial_{x} A^{\eps_{j}} 
+ \partial_{x}^{\eps_{j}} u_{1}^{\eps_{j}}
       \big) + A^{\eps_{j}} \partial_{x}u^{\eps_{j}}_{1} 
+ \frac{ 1}{\eps_{j}^2} g(\eps_{j}^2 A^{\eps_{j}})
        \partial_{x} A^{\eps_{j}}\Big)       \partial_{x} \zeta \\
      & & - \int_{[0, T] \times \mathbb{R}^n}\Delta_{\perp}u_{1}^{\eps_{j}}\zeta 
+ \frac{ 1 }{ 4 c^2} \int_{[0, T] \times \mathbb{R}^n} \partial_{xx} A^{\eps_{j}}  \partial_{xx} \zeta= \int_{\mathbb{R}^n} \partial_{x}\big( A_{0}^{\eps_{j}}
        + (u_{0})_{1}^{\eps_{j}}\big) \zeta(0, X)\, dX + R^{\eps_{j}}
     \end{eqnarray*}
           for every $\zeta \in \mathcal{C}^\infty_{c}(\mathbb{R}\times \mathbb{R}^n)$, where
       thanks to the uniform bound \eqref{Au}, we have
      $$| R^{\eps_{j}}| \leq C \eps_{j}.$$
    We can easily pass to the limit
       in the above formulation by using that in the nonlinear terms one converges
        strongly and one weakly.
       We thus get 
       by using again an expansion of  $g(\eps^2 A^\eps)$, that  
       \begin{eqnarray*}
    & & \int_{[0, T] \times \mathbb{R}^n}\Big(
    2   \partial_{x}A
      \partial_{t} \zeta +  k A \partial_{x} A \partial_{x}\zeta 
       - \Delta_{\perp} A \,  \zeta  +  \frac{ 1}{4 c^2} 
       \partial_{xx} A  \partial_{xx} \zeta\Big) dt dX= \int_{\mathbb{R}^n} \partial_{x}\big( A_{0}
        + (u_{0})_{1}\big) \zeta(0, X)\, dX 
      \end{eqnarray*}
      which is
 the weak form of the KP-I equation (or KdV)
 $$ \partial_{x} \Big( 2   \partial_{t} A  +   k A \partial_{x} A - 
\frac{ 1}{4 c^2 } \partial_{x}^3 A \Big) + \Delta_{\perp} A = 0$$
  with initial value
  $$ A_{|t=0} = \frac{ 1}{2 }\big( A_{0}+  \frac{ 1}{ 2 c }\partial_{x} \varphi_{0}\big).$$
  
  Furthermore,   thanks to the  uniqueness  of $H^s$ solutions, $s> 1+ n/2$ for the KP-I equation, 
   we get that the full sequence $A^\eps$, $\partial_{x} \varphi^\eps$ converges. 
 
 Note that in dimension $1$, we  can get compactness in time by writting directly that
 $$
 \left\{  \begin{array}{ll} \displaystyle{ \partial_{t}  A^\eps  
+ \frac{ 1}{ \eps^2} \partial_{x}\big(  u^\eps_{1} -  A^\eps \big) = S_{A}^\eps}, \\  \\
    \displaystyle{ \partial_{t}u_{1}^\eps + \frac{ 1 }{\eps^2} \partial_{x} \big(A^\eps -u^\eps_{1}
      \big) = S_{u}^\eps} \end{array}\right.$$
with 
  $$
    \big|\!\big|(S_{A}^\eps , S_{u}^\eps)  \big|\!\big|_{H^{-1}} \leq C, \quad \forall t  \in[0, T], 
\quad \forall \eps \in (0, \eps_{0}]  $$
     for some $C>0$  since the  apparently singular term $\eps^{-1} \nabla_{\perp}\cdot  u^\eps_{\perp}$
      is absent in dimension $1$. Then we can finish as in the proof of  Theorem \ref{KdVH1}.
      Thus we get in particular that $A^\eps$ converges strongly towards $A$ in 
      $L^2(0, T, H^{m+1}_{loc})$ (for $n\geq2$
      we have only proven the strong convergence  in $L^2(0, T, H^{m}_{loc})$ for $\partial_{x} A^\eps$).\\

 In the general $n$-dimensional case, it remains to show that, if $ u^\e_\perp=
  \eps \nabla_\perp \varphi^\eps  \to 0$ in $L^2$, then
$$ \frac{1}{2} \big( A^\e + u^\e_1 \big) \to A \quad \quad {\rm in} \quad L^2\big( [0,T], L^2 \big). $$
Indeed, the convergences in $L^2\big( [0,T], H^\s \big)$ for $0\leq \s < s$ will then follow by 
interpolation on space using the bounds \eqref{unif}.

We recall that the scaled energy writes
$$ E^\e (\psi^\e) = \frac{1}{2} \int_{\R^n} 
|\p_x \psi^\e|^2 + \e^2 | \nabla_\perp \psi^\e|^2 
+ \frac{1}{\e^2} F\big( |\psi^\e|^2 \big) \ dX, $$
and we recall the expansion to second order
$$ F (R) = c^2 \big( R -1 \big)^2 + F_3(R ), \quad \quad {\rm with} \quad 
F_3(1 + r ) = \BO(r^3), \quad r\to 0.$$
Moreover, we have, on $[0,T]$,
$$ \psi^\e = \rho^\eps  \exp \big( i\e \vp^\e \big), \quad \rho^\eps = 1+ \e^2 A^\e, $$
and using that  for $1 \leq j\leq n$, $|\p_j \psi|^2 = \e^4 (\p_j A^\e )^2 
+ \e^2 (\rho^\eps)^2 (\p_j \vp^\e)^2$, we infer as in the proof of Lemma \ref{Bornes} 
the following equality:
\begin{align}
\label{devel}
E^\e(\psi^\e) = & \ \frac{\e^2}{2} \int_{\R^n} (\p_x \vp^\e)^2 
+ \frac{c^2}{\e^4} \big( (\rho^\e)^2 - 1 \big)^2 
+ \big( (\rho^\e)^2 - 1 \big) \cdot (\p_x \vp^\e)^2 + \e^2 (\p_x A^\e)^2 \ dX \nonumber \\ 
& + \frac{1}{2} \int_{\R^n} \e^4 (\rho^\e)^2 |\nabla_\perp \vp^\e|^2 
+ \e^2 |\nabla_\perp \rho^\e |^2 + \frac{1}{\e^2}\, F_3 \big( (\rho^\e)^2 - 1 \big) \ dX \\ 
\label{Enerv} = & \ \frac{\e^2}{2} \int_{\R^n} 
(\p_x \vp^\e)^2 + 4c^2 (A^\e)^2 + \e^2 |\nabla_\perp \vp^\e |^2\ dX  + \BO(\e^4)
\end{align}
uniformly on $[0,T]$. To get the last line, we have used   \eqref{unif}, which yields that $|\!| A^\e |\!|_{L^\ii} \leq K $, hence 
$|\!| (\rho^\e)^2 -1 |\!|_{L^\ii} \leq K\e^2$, 
$$ \Big| \int_{\R^n} \big( (\rho^\e)^2 -1 \big) (\p_x \vp^\e)^2 \ dX \Big| \leq K\e^2 
\quad \quad {\rm and} \quad \quad 
\Big| \int_{\R^n} \frac{1}{\e^2}\, F_3\big( (\rho^\e)^2 -1 \big) \ dx \Big| \leq K \e^4.$$

Furthermore, we may define (if $n\geq 2$) the momentum in the $x$ direction by
$$ P^\e \big( \psi^\e \big) \equiv \frac{\e}{2} \int_{\R^n} \big( (\rho^\e)^2 - 1 \big) \p_x \vp^\e \ dX $$
for maps $\psi^\e = \rho^\e e^{i\e \vp^\e }$ with $\rho^\e = | \psi^\e | \geq 1/2$. In view of 
the bounds \eqref{unif}, $| \psi^\e | \geq 1/2$ on $[0,T]$ (for $0< \e \leq \e_0$), hence $\psi^\e$ 
has a well-defined momentum, which is independent of $t\in [0,T]$. Morever, there holds, uniformly on 
$[0,T]$,
\begin{align}
\label{develomoment}
P^\e(\psi^\e) = \frac{\e}{2} \int_{\R^n} \big( (\rho^\e)^2 - 1 \big) \p_x \vp^\e \ dX 
= & \ \frac{\e^2}{2} \int_{\R^n} \big( 2 A^\e + \e^2 (A^\e)^2 \big) \p_x \vp^\e \ dX \nonumber \\
= & \ \e^2 \int_{\R^n} A^\e \p_x \vp^\e\ dX + \BO(\e^2).
\end{align}
As a consequence, in view of \eqref{unif},
$$ E^\e(\psi^\e) + 2c P^\e(\psi^\e) = 2c^2 \e^2 \int_{\R^n} 
\big( A^\e + u^\e_1 \big)^2 +  |u^\e_\perp |^2\ dX + \BO(\e^4)$$
uniformly on $[0,T]$. At the initial time $t=0$, we have
$$ E^\e(\psi^\e_0) + 2c P^\e(\psi^\e_0) = 
2c^2 \e^2 \int_{\R^n} \big( A^\e_0 + (u^\e_0)_1 \big)^2 + | (u^\e_0)_\perp |^2\ dX + \BO(\e^4) ,$$
hence, by conservation of $E^\e(\psi^\e) + 2c P^\e(\psi^\e)$ for $0\leq t \leq T$,
\be
\label{schpoutz}
\int_{\R^n} \big( A^\e(t) + u^\e_1(t) \big)^2 +   \big|u^\e_\perp (t)\big|^2\ dX 
= \int_{\R^n} \big( A^\e_0 + (u^\e_0)_1 \big)^2 +  \big| (u^\e_0)_\perp \big|^2 \ dX + \BO(\e^2),
\ee
uniformly for $t\in [0,T]$. We consider now
$$ \nu^\e \equiv \int_0^T \big|\!\big| A^\e + u^\e_1 - 2A \big|\!\big|^2_{L^2} + 
 \big|\!\big| u^\e_\perp \big|\!\big|^2_{L^2} \ dt .$$
Expansion gives
$$ \nu^\e = \int_0^T \big|\!\big| A^\e + u^\e_1 \big|\!\big|^2_{L^2} +  \big|\!\big| u^\e_\perp \big|\!\big|^2_{L^2} - 4 \big|\!\big| A \big|\!\big|^2_{L^2}  \ dt 
- 4 \int_0^T \langle A^\e + u^\e_1 -2 A, A \rangle_{L^2} \ dt.$$
One can show exactly as in the end of subsect. \ref{finpreuve} that since 
$A\in \BC\big( [0,T], L^2 \big)$ and $A^\e$, $u^\e_1$ converge to $A$ weakly in 
$L^2 ( [0,T],L^2_{loc})$, then
$$ \int_0^T \langle A^\e + u^\e_1 -2 A, A \rangle_{L^2} \ dt \to 0 \quad \quad 
{\rm as} \quad \quad \e \to 0.$$
Moreover, since the $L^2$ norm of the solution $A$ of KP-I does not 
depend on time,
$$ \big|\!\big| 2A(t) \big|\!\big|_{L^2} = \big|\!\big| 2A_{|t=0} \big|\!\big|_{L^2} 
= \big|\!\big| A_0 + (u_0)_1 \big|\!\big|_{L^2}.$$
Hence,  by using  \eqref{schpoutz}, we find after an integration in time that  
$$ \nu^\e = T \Big( \big|\!\big| A^\e_0 + (u^\e_0)_1 \big|\!\big|^2_{L^2} 
- \big|\!\big| A_0 + (u_0)_1 \big|\!\big|^2_{L^2} 
+  \big|\!\big| (u^\e_0)_\perp \big|\!\big|^2_{L^2} \Big) + o(1). $$
Thanks to our assumption \eqref{CVinit}, we thus get  $\nu^\e \to 0$ as required.

 \subsection{Proof of Theorem \ref{H1perp}}
 
 \label{section3}

To use the assumption \eqref{assezprep} in order to get the convergence in stronger norms,  
we will follow the lines of the proof of Lemma \ref{Bornes}. From 
\eqref{devel}, we infer
\begin{align}
\label{develE-cP}
E^\e(\psi^\e) - 2c P^\e(\psi^\e) = & \ \frac{\e^2}{2} \int_{\R^n} 
\big( (\rho^\e)^2-1 \big) (\p_x \vp^\e)^2 + \e^2(\p_x A^\e )^2  \nonumber
\\ & + \Big( \p_x \vp^\e - \frac{c}{\e^2} \big( (\rho^\e)^2-1 \big) \Big)^2
+ \e^2 (\rho^\e)^2 |\nabla_\perp \vp^\e|^2 \ dX  \\
& + \frac{1}{2} \int_{\R^n} \e^6 |\nabla_\perp A^\e|^2
+ \frac{1}{\e^2}\, F_3 \big( (\rho^\e)^2-1 \big) \ dX.\nonumber
\end{align}

Let
$$ \d^\e \equiv \big|\!\big| \p_x \vp^\e_0 - 2cA^\e_0 \big|\!\big|_{L^2} $$
 which tends to zero by assumption.
 As in the proof of  Lemma \ref{Bornes}, we have thanks to \eqref{assezprep}
  in the case $n\geq 2$  the following upper bounds
\be
\label{majoE}
E^\e(\psi^\e_0) = 
\frac{\e^2}{2} \int_{\R^n} 4 c^2 (A_0^\e)^2 + (\p_x \vp_0^\e)^2 \ dX + \BO( \e^4) 
= 4 c^2 \e^2 \int_{\R^n } A_0^2 \ dX + o(\e^2) \leq K\e^2 
\ee
and
$$ E^\e(\psi^\e_0) - 2c P^\e(\psi^\e_0) \leq K \e^4 + \e^2 (\d^\e)^2.$$
Note that here, we have used that 
$$ \Big|\!\Big| \p_x \vp^\e_0 - \frac{c}{\e^2} \big( (\rho^\e_0)^2 -1 \big) \Big|\!\Big|_{L^2} 
= \big|\!\big| \p_x \vp^\e_0 - 2c A^\e_0 - c\e^2 (A^\e_0)^2 \big|\!\big|_{L^2} 
\leq \big|\!\big| \p_x \vp^\e_0 - 2c A^\e_0 \big|\!\big|_{L^2} + 
c\e^2 \big|\!\big| (A^\e_0)^2 \big|\!\big|_{L^2} \leq \d^\e + K \e^2. $$

As a consequence, since $E^\e(\psi^\e)$ and $P^\e(\psi^\e)$ do not depend on time,
\begin{align}
\nonumber
K \e^4 + \e^2 (\d^\e)^2 \geq & \ E^\e \big( \psi^\e(t) \big) - 2c P^\e\big( \psi^\e(t) \big) \\ 
\geq & \ \frac{\e^4}{2} \int_{\R^n} (\p_x A^\e)^2 + (\rho^\e)^2 |\nabla_\perp \vp^\e|^2\ dX 
+ \frac{\e^2}{2} \int_{\R^n} \Big( \p_x \vp^\e -\frac{c}{\e^2} \big( (\rho^\e)^2-1 \big) \Big)^2 
\ dX \\ 
\nonumber& \quad - 
\frac{1}{2} \Big| \int_{\R^n} \big( (\rho^\e)^2 -1 \big) (\p_x \vp^\e)^2 \ dX \Big| - 
\Big| \int_{\R^n} \frac{1}{2\e^2}\, F_3\big( \rho^2-1 \big) \ dX \Big|\\
\label{derniere}\geq & \ \frac{\e^4}{2} \int_{\R^n} (\rho^\e)^2 |\nabla_\perp \vp^\e|^2\ dX
+ \frac{\e^2}{2} \int_{\R^n} \Big( \p_x \vp^\e -\frac{c}{\e^2} \big( (\rho^\e)^2-1 \big) \Big)^2 
\ dX - K \e^4.
\end{align}
This gives the estimate
\be
\label{smoothA-u} \sup_{0 \leq t \leq T} 
\int_{\R^n} \Big( \p_x \vp^\e -\frac{c}{\e^2} \big( (\rho^\e)^2-1 \big) \Big)^2 \ dX
\leq K \e^2 + 2 (\d^\e)^2 \to 0 \quad \quad {\rm as} \quad \e \to 0\ee
 in all dimensions $n\geq 1$.

Furthermore, in dimension $n\geq 2$, since  $\d^\e =\BO(\e)$,   we also get from \eqref{derniere} 
that 
$$ \int_{\R^n} (\rho^\e)^2 |\nabla_\perp \vp^\e|^2\ dX \leq K.$$
Thus, we have obtained \eqref{gradperp} since $\rho^\e \geq 1/2$.  \\
\bigskip

 From \iref{euler1}, $A^\e + u^\e_1$ solves
$$ \p_t \big( A^\e + u^\e_1 \big) + 2 u^\e \cdot \nabla^\e \big( A^\e + u^\e_1 \big) 
+ (k-5) A^\e \p_x A^\e + A^\e \nabla^\e \cdot u^\e + \Delta_\perp \vp^\e
= \p_x \Big( \frac{\Delta^\e A^\e}{4 c^2 \rho^\e} \Big). $$
In view of the  the $H^s$ bounds \eqref{unif} in Theorem \ref{Smoothy}, and possibly \iref{gradperp} 
if $n\geq 2$, we then infer
\be
\label{Kompactus}
\big|\!\big| A^\e + u^\e_1 \big|\!\big|_{\BC ( [0,T], H^s )} \leq K 
\quad \quad \quad {\rm and} \quad \quad \quad 
\big|\!\big| \p_t \big( A^\e + u^\e_1 \big) \big|\!\big|_{L^\ii ( [0,T], H^{-1} )} \leq K.
\ee
This implies, by Aubin-Lions's Lemma (see, {\it e.g.}, \cite{S}), that for any 
$0 \leq \s < s$, $A^\e+ u^\e_1$ is precompact in $\BC \big( [0,T], H^\s_{loc} \big)$. 
 From \eqref{smoothA-u}, we know that 
$$ \p_x \vp^\e - 2 c A^\e = 2 c \big(  u^\e_1 - A^\e \big)\to 0 
\quad \quad {\rm in} \quad \BC \big( [0,T] , L^2 \big). $$
Combining this with the $H^s$ bounds \eqref{unif}, this yields, 
by interpolation, for $0\leq \s < s$,
$$ \p_x \vp^\e - 2 c A^\e = 2 c \big(  u^\e_1 - A^\e \big)\to 0 
\quad \quad {\rm in} \quad \BC \big( [0,T] , H^\sigma \big). $$
In particular,
$$ A^{\e} \to A \quad \quad \quad {\rm and} \quad \quad \quad \p_x \vp^{\e} \to 2c A 
\quad \quad {\rm in} \quad \BC \big( [0,T], H^\sigma_{loc} \big). $$
We can now prove that, as $\e \to 0$,
$$ A^\e \to A  \quad \quad {\rm in} \quad \BC \big( [0,T], L^2 \big). $$
Indeed, we may follow the lines of  the end of the proof of Theorem \ref{KdVH1} in 
Sect. \ref{finpreuve} since thanks to \eqref{unif},  \eqref{gradperp} (if $n=2,3$) and 
\eqref{smoothA-u}, the expansion
$$ E^\e(\psi^\e) = \frac{\e^2}{2} \int_{\R^n} 4c^2 \big( A^\e(t) \big)^2 + \big( \p_x \vp^\e(t) \big)^2 \ dX 
+ \BO(\e^4) = 4c^2 \int_{\R^n} \big( A^\e(t) \big)^2 \ dX + o(\e^2) $$
holds  uniformly for $0 \leq t \leq T$  and 
$\ds{\BI_0\big( A(t) \big) = |\!| A(t) |\!|^2_{L^2} } = |\!| A_0 |\!|^2_{L^2}$ do not depend on 
$t\in [0,T]$. Notice indeed that in this case, the initial datum for KP-I is
$$ A_{|t=0} = \frac{1}{2} \big( A_0 + \frac{1}{2c} \p_x \vp_0 \big) = A_0. $$

From the $H^s$ bounds  \eqref{unif}  and by interpolation in space, we  finally get that 
$$ \forall \ 0 \leq \s < s \quad \quad A^\e \to A \quad \quad    {\rm  in} \quad \BC \big( [0,T], H^{\sigma+ 1} \big)
 \quad \quad {\rm and} \quad \quad \p_x \vp^\e \to 2c  A
\quad \quad {\rm in} \quad \BC \big( [0,T], H^\s \big).$$




\begin{thebibliography}{99}


\bibitem{Alazard-Carles}{\sc T. Alazard and R. Carles}, 
{\it Supercritical geometric optics for Nonlinear Schrodinger equations}, Preprint 2007.


\bibitem{AL}{\sc B. Alvarez-Samaniego and D. Lannes},
{\it Large time existence for 3D water-waves and asymptotics.}
Invent. Math. {\sc 171}(2008), no. 3, 485-541.

\bibitem{BR} {\sc N. Berloff and P. Roberts,} 
{\it Motions in a Bose condensate: X. New results on stability of axisymmetric solitary waves of the Gross-Pitaevskii equation.} 
J. Phys. A: Math. Gen., {\bf 37} (2004), 11333-11351.

\bibitem{BDS} {\sc F. B\'ethuel, R. Danchin and D. Smets,} 
{\it On the linear wave regime of the Gross-Pitaevskii equation.} Preprint.

\bibitem{BGS} {\sc F. B\'ethuel, P. Gravejat and J-C. Saut,} 
{\it On the KP I transonic limit of two-dimensional Gross-Pitaevskii travelling waves.} 
Dynamics of PDE {\bf 5}, 3 (2008), 241-280.

\bibitem{BGSS} {\sc F. B\'ethuel, P. Gravejat, J-C. Saut and D. Smets,} 
{\it On the Korteweg-de Vries long-wave transonic approximation
of the Gross-Pitaevskii equation.} Preprint.

\bibitem{Chiron-Rousset}{\sc D. Chiron and F. Rousset},
{\it Geometric optics and boundary layers for Nonlinear Schrodinger equations}, Preprint 2008.

\bibitem{Ga} {\sc C. Gallo,} 
{\it The Cauchy Problem for defocusing Nonlinear Schrödinger equations with 
non-vanishing initial data at infinity.} 
Comm. Partial Differential Equations {\bf 33}, no. 4-6 (2008), 729-771.

\bibitem{Ge} {\sc P. G\'erard,} 
{\it The Gross-Pitaevskii equation in the energy space.} 
Stationary and Time Dependent Gross-Pitaevskii Equations", A. Farina and 
J.-C. Saut editors, Contemporary Mathematics, American Mathematical Society (2008).

\bibitem{GP} {\sc V. Ginzburg and L. Pitaevskii,} 
{\it On the theory of superfluidity.} 
Sov. Phys. JETP {\bf 34} (1958), 1240.

\bibitem{Grenier}{\sc E. Grenier,} 
{\it Semiclassical limit of the nonlinear Schr\"odinger equation in small time.} 
Proc. Amer. Math. Soc. {\bf 126} (1998), no. 2, 523--530.

\bibitem{GrenierS} {\sc E. Grenier,}
{\it Pseudo-differential energy estimates of singular perturbations.}
 Comm. Pure Appl. Math. {\bf 50} (1997), no. 9, 821--865.

\bibitem{G} {\sc E. Gross,} 
{\it Hydrodynamics of a superfluid condensate,} 
J. Math. Phys. {\bf 4}, (2) (1963), 195-207.

\bibitem{IKT} {\sc A. Ionescu, C. Kenig and D. Tataru,} 
{\it Global well-posedness of the KP-I initial-value problem in the energy space.} 
Invent. Math. {\bf 173}, no. 2 (2008), 265-304.

\bibitem{JR} {\sc C. Jones and P. Roberts,}
{\it Motion in a Bose condensate: IV. Axisymmetric solitary waves.}
J. Phys. A: Math. Gen., {\bf 15} (1982) 2599-2619.

\bibitem{KPV} {\sc C. Kenig, G. Ponce and L. Vega,} 
{\it Well-posedness of the initial value problem for the Korteweg-de Vries equation.} 
J. Amer. Math. Soc. {\bf 4} (1991), no. 2, 323-347.

\bibitem{KPV2} {\sc C. Kenig, G. Ponce and L. Vega,} 
{\it A bilinear estimate with applications to the KdV equation} 
J. Amer. Math. Soc. {\bf 9} (1996), no. 2, 573-603.

\bibitem{KAL} {\sc Y. Kivshar, D. Anderson and M. Lisak,}
{\it Modulational instabilities and dark solitons in a generalized nonlinear Schr\"odinger-equation.} 
Phys. Scr. {\bf 47}, (1993) 679-681.

\bibitem{KL} {\sc Y. S. Kivshar and B. Luther-Davies}, 
{\it Dark optical solitons: physics and applications}. 
Physics Reports {\bf 298} (1998), 81-197.

\bibitem{Klainerman-Majda} {\sc  S. Klainerman and A. Majda}
{\it Singular limits of quasilinear hyperbolic systems with large parameters and the incompressible limit of compressible fluids.} Comm. Pure Appl. Math. {\bf 34} (1981), no. 4, 481--524.

\bibitem{L} {\sc Y. Liu,} 
{\it Strong instability of solitary-wave solutions to a Kadomtsev-Petviashvili 
equation in three dimensions.} 
J. Differential Equations, {\bf 180} no. 1 (2002), 153-170.

\bibitem{MST} {\sc L. Molinet, J.-C. Saut and N. Tzvetkov,} 
{\it Global well-posedness for the KP-I equation.} 
Math. Ann. {\bf 324}, no. 2 (2002), 255-275.

\bibitem{RB} {\sc P. Roberts and N. Berloff,} 
{\it Nonlinear Schr\"odinger equation as a model of superfluid helium.} 
In "Quantized Vortex Dynamics and Superfluid Turbulence" edited by C.F. Barenghi, 
R.J. Donnelly and W.F. Vinen, Lecture Notes in Physics, volume 571, Springer-Verlag, 2001.

\bibitem{Schochet}{\sc S.  Schochet,} 
{\it  Asymptotics for symmetric hyperbolic systems with a large parameter.}
 J. Differential Equations {\bf 75} (1988), no. 1, 1--27.

\bibitem{S} {\sc J. Simon,} 
{\it Compact sets in the space $L^p(0,T;B)$.} 
Ann. Mat. Pura Appl. {\bf 146}, (4) (1987), 65-96.

\bibitem{Tsu} {\sc T. Tsuzuki,}
{\it Nonlinear waves in the Pitaevskii-Gross equation,}
J. Low Temp. Phys. {\bf 4}, no. 4 (1971) 441-457.

\bibitem{Z} {\sc P. Zhidkov,} {\it Korteweg-De Vries and nonlinear 
Schr\"odinger equations : qualitative theory,} Volume 1756 of Lecture Notes in Mathematics. Springer-Verlag, Berlin, 2001.

\bibitem{Zhou}{\sc Y. Zhou,}
 {\it Uniqueness of weak solution of the KdV equation}. 
Internat. Math. Res. Notices (1997), no. 6, 271--283. 

\end{thebibliography}
\end{document}